\RequirePackage{iftex}
\ifXeTeX
 \PassOptionsToPackage{preserveurlmacro}{breakurl}
 \makeatletter
 \providecommand{\headerps@out}[1]{}
 \makeatother
\fi
\PassOptionsToPackage{noabbrev}{cleveref}
\documentclass[hidelinks,nohypdvips]{siamart251216}

\usepackage[T1]{fontenc}
\usepackage{amsmath,amssymb,mathtools}
\usepackage{enumitem}
\usepackage{microtype}
\usepackage{algorithm}
\usepackage{algpseudocode}
\usepackage{booktabs}
\usepackage{siunitx}
\usepackage{needspace}

\newsiamremark{remark}{Remark}
\newsiamthm{assumption}{Assumption}

\crefname{assumption}{Assumption}{Assumptions}
\AddToHook{env/assumption/begin}{\crefalias{theorem}{assumption}}
\AddToHook{env/remark/begin}{\crefalias{theorem}{remark}}

\newcommand{\A}{\mathcal{A}}
\newcommand{\Aad}{\mathcal{C}_{\mathrm{ad}}}
\newcommand{\Aadomega}{\mathcal{C}_{\mathrm{ad},\omega}}
\newcommand{\R}{\mathbb{R}}
\newcommand{\dd}{\,\mathrm{d}}
\newcommand{\cl}{\mathrm{cl}}
\newcommand{\diag}{\operatorname{diag}}
\newcommand{\Dom}{\operatorname{dom}}

\newcommand{\Ls}{\mathcal{L}}
\newcommand{\norm}[1]{\left\|#1\right\|}
\newcommand{\ip}[2]{\left(#1,#2\right)}
\newcommand{\duality}[2]{\left\langle #1,#2\right\rangle}
\newcommand{\KA}{\mathcal{K}_A}
\newcommand{\KB}{\mathcal{K}_B}

\newcommand{\LA}{\mathfrak L_A}

\newcommand{\revise}[2]{{\color{red}#1}\ifmmode\,\else\space\fi{\color{blue}#2}}

\headers{Optimal Actuator Support Design}{B. Azmi and D. Wachsmuth}

\title{Optimal Actuator Support Design for Infinite-Horizon Linear-Quadratic Parabolic Control Problems}
\author{BEHZAD AZMI\thanks{Department of Mathematics and Statistics, University of Konstanz, 78457 Konstanz, Germany. Email: \texttt{behzad.azmi@uni-konstanz.de}} \and
DANIEL  WACHSMUTH \thanks{Institute of Mathematics, University of Würzburg, 97074 Würzburg, Germany, Email: \texttt{daniel.wachsmuth@uni-wuerzburg.de}. The research of the second author was partially funded by the Deutsche Forschungsgemeinschaft (DFG, German Research Foundation) -- project
number 565876882.}}
\date{\today}

\begin{document}

\maketitle

\begin{abstract}
We study an optimal actuator design problem for infinite-horizon LQ control
problems governed by parabolic equations. The control acts only on a small
measurable subregion of the spatial domain, and the objective combines the
infinite-horizon control value with a cost for the actuator region. Using the Riccati representation of the value function, we derive a comparison identity for two actuator supports and obtain an exact expansion of the objective. This leads to a topological derivative with respect to local changes in the support. We prove a localized Riccati sensitivity estimate to control the remainder and show that the infinite-horizon term in the derivative can be computed by solving a Lyapunov equation. This yields a computable support-update algorithm with a measure-controlled Armijo search. Under explicit uniform remainder and admissibility assumptions, its first-order stationarity defect tends to zero. Numerical experiments illustrate these results.
\end{abstract}


\section{Introduction}
\label{sec:introduction}

The goal of this paper is to study an actuator-design problem for
infinite-horizon linear-quadratic control problems governed by parabolic
equations. The central question is not only how to compute an optimal control
for a fixed actuator configuration, but also how to choose the region in the
spatial domain where the control is allowed to act. This leads to an
optimization problem with respect to the actuator support itself.

Let \(\Omega\subset\R^d\), with $d\leq 3$, be a bounded domain and let \(A\subset\Omega\)
denote the actuator region. A typical model is the controlled
convection--reaction--diffusion equation
\begin{equation}
\label{eq:motivating_pde}
\begin{gathered}
    \partial_t y-\kappa\Delta y+b\cdot\nabla y+cy
    =\chi_Au,
    \qquad t>0,\ x\in\Omega,\\
    \mathcal G y=0\quad\text{on }\partial\Omega,
    \qquad y(0)=y_0\quad\text{in }\Omega.
\end{gathered}
\end{equation}
Here \(y=y(t,x)\) is the state, \(u=u(t,x)\) is the control,
\(\kappa>0\) is the diffusion coefficient, \(b=b(x)\) is a convection field,
and \(c=c(x)\) is a reaction coefficient. We choose
\(\mathcal G y=y|_{\partial\Omega}\) or
\(\mathcal G y=\partial_\nu y|_{\partial\Omega}\) for homogeneous Dirichlet
or Neumann boundary conditions, respectively, with \(\nu\) the outward unit
normal. This choice is fixed throughout the support optimization.
The characteristic function \(\chi_A\) restricts the control action to \(A\).

For a fixed support \(A\), the control may be physically seen as an element
of \(L^2(A)\), extended by zero to the entire domain.

However, since the support \(A\) is itself the design variable, it is more convenient for the analysis to work on the fixed ambient control space $U:=L^2(\Omega)$ and to represent the actuator through the multiplication operator $ u\mapsto \chi_Au$. This avoids changing the control space when the support changes.

For a fixed actuator support \(A\), the associated infinite-horizon
linear-quadratic cost is
\begin{equation}
\label{eq:motivating_cost}
    \mathcal J_A(u;y_0)
    :=
    \int_0^\infty
    \left(
        \norm{Qy(t)}_{H_Q}^2
        +
        \alpha\norm{u(t)}_{L^2(\Omega)}^2
    \right)\dd t,
\end{equation}
where \(y\) solves \eqref{eq:motivating_pde}, \(Q\) is an observation or
state-cost operator, and \(\alpha>0\) is the control-regularization parameter.
The corresponding value function is defined by
\begin{equation}
\label{eq:value_function}
    V_A(y_0)
    :=
    \inf_{u\in L^2(0,\infty;L^2(\Omega))}
    \mathcal J_A(u;y_0),
\end{equation}
subject to \eqref{eq:motivating_pde}.
% Thus,  \(V_A(y_0)\) is the optimal infinite-horizon LQ value associated with the actuator support \(A\) and the initial condition \(y_0\).
The actuator-design problem then consists in minimizing over both the control and the support \(A\). In this paper, we consider objectives of the form
\begin{equation}
\label{eq:motivating_support_cost}
    J(A)
    =
    V_A(y_0)
    +
    \int_\Omega \chi_A(x)\beta(x)\dd x,
\end{equation}
where the second term penalizes the actuator's size or spatial cost. The
weight \(\beta\) can be used, for instance, to encode regions where actuation
is expensive or undesirable.

The main analytical challenge is that changing \(A\) changes the control
operator, the Riccati operator, the optimal feedback, and the closed-loop
semigroup. We aim to derive a first-order expansion of the
support-dependent objective under local changes of the actuator region. The
leading term of this expansion yields a topological derivative, which can be
used to decide where the support should be enlarged or removed. A key point is
that the infinite-horizon contribution to this derivative can be computed
through a Lyapunov equation, so that the resulting descent method does not
require direct integration over an infinite time interval.

The subsequent analysis is carried out in an abstract Hilbert-space setting.
This allows us to treat not only \eqref{eq:motivating_pde}, but more generally
linear parabolic equations generated by forms satisfying a G{\aa}rding inequality on a Gelfand triple.
Diffusion--convection--reaction equations of the form above are included under
the usual boundedness and ellipticity assumptions.

\paragraph{Related literature and contributions}
Actuator location via Riccati equations for infinite-dimensional linear systems is
studied in \cite{Morris2011}, using worst-case or averaged initial-state
performance. That work also establishes the convergence of optimal costs and
locations under suitable finite-dimensional approximations. For structural
models, \cite{DarivandiMorrisKhajepour2013} develops a branch-and-bound subgradient-based method for optimal LQ placement among candidate locations, utilizing convexity in actuator weights.
We optimize a distributed actuator support as a set-valued design
variable.

For semilinear distributed parameter systems,
\cite{EdalatzadehMorris2019} investigates joint optimization of the control and
the actuator operator on a finite time interval. Under suitable compactness
and differentiability assumptions, it proves the existence of optimal designs
and controls and derives first-order necessary conditions. This allows
nonlinear dynamics, whereas our analysis concerns autonomous linear systems
over an infinite horizon.

Particularly close to our formulation, \cite{KaliseKunischSturm2018}
studies parabolic actuator design for a prescribed initial state and for
worst-case performance, deriving shape and topological sensitivities and
corresponding level-set methods. Its sensitivity analysis is finite-horizon,
while the numerical experiments in \cite[Sec.~6]{KaliseKunischSturm2018}
use infinite-horizon LQR values computed from algebraic Riccati equations.
For the Euler--Bernoulli vibration model,
\cite{EdalatzadehKaliseMorrisSturm2022} combines finite-horizon LQR
performance with shape and topological sensitivities and a level-set
algorithm.

In the present work, we derive the support sensitivity of the infinite-horizon parabolic LQR value directly. The objective also includes the spatial cost \(\int_A\beta(x)\dd x\) of the actuator region. We do not assume that the uncontrolled equation is exponentially stable. Instead, we restrict the optimization to supports \(A\) for which a
feedback control acting in \(A\) makes the closed-loop system
exponentially stable. A Riccati comparison identity and a localized
sensitivity estimate yield the first-order support expansion. Its
infinite-horizon term is evaluated through a closed-loop Lyapunov equation,
without time truncation. Thus, the distinction from \cite{KaliseKunischSturm2018} lies in the direct infinite-horizon sensitivity analysis, not in the use of topology optimization itself.

Our descent method updates the actuator support directly using
measure-controlled Armijo steps, following the support formulation for \(L^0(\Omega)\)-penalized control in \cite{Wachsmuth2024}, rather than evolving a level-set function.
Their application here requires estimates for the support-dependent Riccati operator and closed-loop dynamics. We therefore state the stationarity analysis under explicit admissibility and uniform remainder assumptions: well-posedness for each support alone does not provide the uniform estimates needed along the sequence.

\paragraph{Organization of the paper}
\Cref{sec:problem_riccati} introduces the problem formulation and collects preliminary results on the Riccati equation, support perturbations,
and the associated perturbation estimates. \Cref{sec:riccati_sensitivity} develops the main analytical tools for the support sensitivity analysis,
including Riccati smoothing, local preservation of stabilizability, and localized sensitivity expansions. \Cref{sec:topological_algorithm} derives the
topological derivative, introduces the corresponding descent algorithm, and analyzes its convergence and stationarity properties.
Finally, \cref{sec:numerics} presents numerical experiments illustrating the algorithm and the theoretical results.

\section{Problem formulation and Riccati preliminaries} \label{sec:problem_riccati}

\subsection{Abstract problem formulation} \label{subsec:abstract_problem}
We set \(H:=L^2(\Omega)\), where \(\Omega\subset\R^d\) is a bounded
domain. We work with the Gelfand triple
\(V\hookrightarrow H\simeq H'\hookrightarrow V'\),
where \(V\) is a Hilbert space and the embeddings are continuous and dense.
A typical example is \(V=H_0^1(\Omega)\).

For \(T\in(0,\infty]\), we set
\(W(0,T):=\{v\in L^2(0,T;V):\partial_t v\in L^2(0,T;V')\}\),
where \(\partial_t v\) denotes the weak time derivative.
Then \(W(0,T)\hookrightarrow C([0,T];H)\) for finite \(T\), and
\(W(0,\infty)\hookrightarrow C([0,\infty);H)\). For \(v\in W(0,T)\), the chain rule
\(\frac{\dd}{\dd t}\norm{v(t)}_H^2
=2\duality{\partial_t v(t)}{v(t)}_{V',V}\) holds almost everywhere;
see e.g.,~\cite{LionsMagenes1972}.

Let \(a:V\times V\to\R\) be a continuous bilinear form satisfying the
G{\aa}rding inequality. Thus, there exist constants \(c_a,\eta>0\) and
\(\lambda\geq0\) such that
\begin{equation}
\label{eq:form_bounds}
    |a(v,w)|
    \leq c_a\norm{v}_V\norm{w}_V,
    \qquad
    a(v,v)
    \geq \eta\norm{v}_V^2-\lambda\norm{v}_H^2,
    \qquad \forall v,w\in V.
\end{equation}
We define \(\A\in\Ls(V,V')\) by
\(\duality{\A v}{w}_{V',V}=-a(v,w)\) for all \(v,w\in V\).
We assume that its realization in \(H\) is closed and densely defined. Arguing as in \cite[Thm.~7.2.7]{Pazy1983} one finds that by \eqref{eq:form_bounds}
the operator $\A$ generates an analytic semigroup on \(H\).

We use \(\A\) for both the variational operator in \(\Ls(V,V')\) and
its part \(\A:\Dom(\A)\subset H\to H\) in \(H\), whenever no
confusion can arise. Then,  for every measurable support \(A\subset\Omega\), we consider the localized
control system
\begin{equation}
\label{eq:state}
    \dot y(t)
    =
    \A y(t)+\chi_Au(t),
    \qquad
    y(0)=y_0.
\end{equation}
Here \(y_0\in H\), \(u\in L^2(0,\infty;H)\), and \(\chi_A\) denotes the characteristic function of \(A\). Thus the control acts only on the
support \(A\). We use the cost \eqref{eq:motivating_cost} and value function
\eqref{eq:value_function}, where \(Q\in\Ls(H,H_Q)\) is bounded and
\(H_Q\) is a Hilbert space. For \(\beta\in L^1(\Omega)\), the actuator-support optimization
problem with objective \eqref{eq:motivating_support_cost} is
\begin{equation}
\label{eq:support_problem}
    \min_{A\in\Aad}J(A),
\end{equation}
where \(\Aad\) is a class of admissible measurable subsets of \(\Omega\).
The admissible class will be chosen so that the infinite-horizon LQR
problem is well-defined for every \(A\in\Aad\).

\subsection{Infinite-horizon Riccati theory and closed-loop regularity} \label{subsec:riccati_theory}

We recall the fixed-support infinite-horizon LQR theory that will be used
throughout the paper.

\begin{assumption}[Stabilizability and detectability]
\label{ass:lqr}
For every \(A\in\Aad\), the pair \((\A,\chi_A)\) is exponentially
stabilizable. The pair \((Q,\A)\) is exponentially detectable.
\end{assumption}

For a fixed support \(A\), the problem is a standard infinite-horizon
linear-quadratic control problem. Under stabilizability and detectability, the
value function is represented by a nonnegative self-adjoint Riccati operator and the
optimal control is a stabilizing feedback; see, e.g.,
\cite{BensoussanDaPratoDelfourMitter2007,CurtainZwart1995,LasieckaTriggiani2000}.

\begin{theorem}
\label{thm:riccati}
Let \cref{ass:lqr} hold. Then, for every \(A\in\Aad\), there exists
a unique nonnegative self-adjoint Riccati operator \(\Pi_A\in\Ls(H)\) such that
\begin{equation}
\label{eq:value_riccati}
    V_A(y_0)
    =
    \ip{\Pi_Ay_0}{y_0}_H
    \qquad
    \text{for all }y_0\in H.
\end{equation}
The optimal control is given in feedback form by
\begin{equation}
\label{eq:feedback_actual}
    u_A^*(t)=\KA y_A(t),
    \qquad
    \KA:=-\frac1\alpha\chi_A\Pi_A.
\end{equation}
The corresponding closed-loop state satisfies
\begin{equation}
\label{eq:closed_loop_A}
    \dot y_A(t)
    =
    \A_A^{\cl}y_A(t),
    \qquad
    y_A(0)=y_0,
\end{equation}
where
\begin{equation}
\label{eq:closed_loop_generator}
    \A_A^{\cl}
    :=
    \A+\chi_A\KA
    =
    \A-\frac1\alpha\chi_A\Pi_A .
\end{equation}
We denote the associated closed-loop semigroup by
\(S_A^{\cl}(t):=e^{t\A_A^{\cl}}\). It is exponentially stable; that is, there
exist constants \(M_A\geq1\) and \(\gamma_A>0\) such that
\begin{equation}
\label{eq:closed_loop_exp_stability}
    \norm{S_A^{\cl}(t)}_{\Ls(H)}
    \leq
    M_Ae^{-\gamma_A t},
    \qquad
    t\geq0.
\end{equation}
Moreover, \(\Pi_A\) satisfies the algebraic Riccati equation in the usual
weak/variational sense. Formally, this equation can be written as
\begin{equation}
\label{eq:riccati_operator}
    \A^*\Pi_A+
    \Pi_A\A-
    \frac1\alpha\Pi_A\chi_A\Pi_A+
    Q^*Q=0.
\end{equation}
\end{theorem}
\begin{proof}
The proof is given in  \cite[Thm.~6.2.7]{CurtainZwart1995}.
\end{proof}

\begin{remark}[Closed-loop Lyapunov form of the Riccati equation]
\label{rem:closed_loop_lyapunov_riccati}
Using \eqref{eq:feedback_actual} and \eqref{eq:closed_loop_generator}, the
Riccati equation \eqref{eq:riccati_operator} can equivalently be rewritten in
closed-loop Lyapunov form as
\begin{equation}
\label{eq:riccati_operator2}
    (\A_A^{\cl})^*\Pi_A
    +
    \Pi_A\A_A^{\cl}
    +
    Q^*Q
    +
    \alpha\KA^*\KA
    =
    0 .
\end{equation}
This identity is the starting point for the closed-loop representation of
\(\Pi_A\) used later in the Riccati smoothing argument.
\end{remark}

\begin{proposition}
\label{prop:parabolic_is_nice}
Let $\Omega\subset\mathbb{R}^d$, $d\leq3$, be a bounded $C^2$-domain,
let $\kappa>0$, $b\in L^\infty(\Omega;\mathbb{R}^d)$,
$c,\beta\in L^\infty(\Omega)$, and $\alpha>0$. Consider $\A=\kappa\Delta-b\cdot\nabla-c$
and let $A\subset\Omega$ be any nonempty open actuator region. Then $(\mathcal{A},\chi_A)$ is exponentially stabilizable. Consequently, for every $y_0\in L^2(\Omega)$, the full-state infinite-horizon LQR problem associated with \eqref{eq:motivating_pde}, and \eqref{eq:motivating_cost} with $H_Q =H=L^2(\Omega )$ and $Q = I$ has finite value and admits a unique optimal
control and a bounded, self-adjoint, nonnegative stabilizing Riccati operator
$\Pi_A\in\Ls(L^2(\Omega))$.
\end{proposition}

\begin{proof}
For the Dirichlet and Neumann realizations of
\(\kappa\Delta-b\cdot\nabla-c\), standard parabolic theory gives a strongly
continuous semigroup \((S(t))_{t\geq0}\) on \(L^2(\Omega)\). The global
Carleman observability theory for linear parabolic equations with bounded
lower-order terms gives, for every nonempty open \(A\subset\Omega\) and every
\(T>0\), the following estimate for the backward adjoint equation. Namely, if
\(\varphi_T\in L^2(\Omega)\) and
\[
    -\frac{\dd}{\dd t}\varphi(t)=\A^*\varphi(t),
    \qquad
    \varphi(T)=\varphi_T,
\]
where \(\A\) is the \(L^2\)-realization of
\(\kappa\Delta-b\cdot\nabla-c\) with the chosen boundary condition, then
\[
    \norm{\varphi(0)}_{L^2(\Omega)}^2
    \leq
    C_T
    \int_0^T
    \norm{\chi_A\varphi(t)}_{L^2(\Omega)}^2\dd t .
\]
See e.g., \cite{FursikovImanuvilov1996,FernandezCaraGuerrero2006} for the Dirichlet
case. For Fourier boundary conditions, which include the Neumann case, see
\cite[Thm.~2]{FernandezCaraGonzalezBurgosGuerreroPuel2006},
which explicitly includes bounded convection and reaction coefficients and
the quantitative null-control estimate; homogeneous Neumann conditions
correspond to a zero boundary Fourier coefficient.
%The diffusion coefficient \(\kappa>0\) is accommodated by rescaling time.
By the Hilbert-space
duality between adjoint observability and null controllability, see e.g., \cite[Sec.~11.2]{TucsnakWeiss2009}, this gives, for every
\(y_0\in L^2(\Omega)\), a control \(u_T\in L^2(0,T;L^2(A))\) such that
\(y(T)=0\) and
\begin{equation}
\label{eq:model_null_control_estimate}
    \norm{u_T}_{L^2(0,T;L^2(A))}
    \leq
    C_T\norm{y_0}_{L^2(\Omega)} .
\end{equation}
Here \(C_T\) may depend on \(T\), \(\Omega\), \(A\), \(\kappa\),
\(\norm{b}_{L^\infty}\), \(\norm{c}_{L^\infty}\), and on the chosen boundary
condition, but it is independent of \(y_0\). We extend $u_T$ by zero from $A\times(0,T)$ to $\Omega \times (0,+\infty)$.

Since the state satisfies \(y(T)=0\), uniqueness for the homogeneous evolution
implies \(y(t)=0\) for all \(t\geq T\).
It follows that $\mathcal J_A(u_T;y_0) < +\infty$, i.e., the
infinite-horizon LQR problem is optimizable. The standard infinite-dimensional
LQR theorem \cite[Thm.~6.2.4]{CurtainZwart1995} therefore gives existence
and uniqueness of the optimal control and a bounded, self-adjoint,
nonnegative operator \(\Pi_A\) representing the value function.

The above null controllability estimate also implies exponential
stabilizability of \((\A,\chi_A)\) by Datko's theorem
\cite[Thm.~V.1.8]{EngelNagel2000};
see also \cite[Ch.~IV, Sec.~4]{Zabczyk1992}.
Since the stage cost is the full
\(L^2(\Omega)\)-norm, the observation operator is the identity, and the
detectability condition is automatic. Hence the stabilizing Riccati theorem
\cite[Thm.~6.2.7]{CurtainZwart1995} gives the stabilizing Riccati operator,
the feedback formula, and exponential closed-loop stability.
\end{proof}

%%%%%%%%%%%%%%%%%%%%%%%%%%%%%%%%%%%%%%%%%%%%%%%%%%%%%%%%%%%%%%%%%%%%%%%%%%%%%

\Cref{thm:riccati} gives the optimal closed-loop trajectory in semigroup
form $y_A(t)=S_A^{\cl}(t)y_0$. We next derive the regularity estimates for \(y_A\) that will be used later.

\begin{lemma}
\label{lem:closed_loop_W_regular}
Let \(A\in\Aad\), and let \(y_A(t)=S_A^{\cl}(t)y_0\) be the optimal
closed-loop state associated with \(y_0\in H\). Then \(y_A\in W(0,\infty)\)
and \(u_A^*=\KA y_A\in L^2(0,\infty;H)\). Moreover, the following estimate
holds:
\begin{equation}
\label{eq:closed_loop_W_estimate}
    \norm{y_A}_{L^2(0,\infty;H)}
    +
    \norm{y_A}_{L^2(0,\infty;V)}
    +
    \norm{\dot y_A}_{L^2(0,\infty;V')}
    +
    \norm{u_A^*}_{L^2(0,\infty;H)}
    \leq
    C_A\norm{y_0}_H,
\end{equation}
where \(C_A>0\) is independent of \(y_0\).
\end{lemma}

\begin{proof} Using the exponential stability estimate \eqref{eq:closed_loop_exp_stability}, we have
\[
    \norm{y_A(t)}_H
    =
    \norm{S_A^{\cl}(t)y_0}_H
    \leq
    M_Ae^{-\gamma_A t}\norm{y_0}_H,
\]
which implies $\norm{y_A}_{L^2(0,\infty;H)}^2\leq  C_A\norm{y_0}_H^2$. Since \(\KA\in\Ls(H)\), it follows immediately that %$u_A^*=\KA y_A\in L^2(0,\infty;H)$, and
$ \norm{u_A^*}_{L^2(0,\infty;H)} \le C_A\norm{y_0}_H$. It remains to show that \(y_A\in L^2(0,\infty;V)\) and
\(\dot y_A\in L^2(0,\infty;V')\). Using \eqref{eq:form_bounds} and testing \eqref{eq:closed_loop_A} with \(y_A\) gives us for a.e.\@ \(t>0\) that
\[
\begin{aligned}
    \frac12\frac{\dd}{\dd t}\norm{y_A(t)}_H^2
    &=
    \duality{\A y_A(t)}{y_A(t)}_{V',V}
    +
    \ip{\KA y_A(t)}{y_A(t)}_H
    \\
    &\leq
    -\eta\norm{y_A(t)}_V^2
    +
    \lambda\norm{y_A(t)}_H^2
    +
    \norm{\KA}_{\Ls(H)}
    \norm{y_A(t)}_H^2.
\end{aligned}
\]
Integrating over \((0,T)\), we obtain
\[
\begin{aligned}
    \eta\int_0^T\norm{y_A(t)}_V^2\,\dd t
    &\leq
    \frac12\norm{y_0}_H^2
    +
    \left(\lambda+\norm{\KA}_{\Ls(H)}\right)
    \int_0^T\norm{y_A(t)}_H^2\,\dd t .
\end{aligned}
\]
Using the already established
\(L^2(0,\infty;H)\)-estimate yields $\norm{y_A}_{L^2(0,\infty;V)} \leq C_A\norm{y_0}_H$.
Finally,
the bound of $\norm{\dot y_A}_{L^2(0,\infty;V')}$ follows from the continuity of \(\A:V\to V'\).
\end{proof}

\subsection{Support perturbations and auxiliary estimates}
\label{subsec:support_perturbations}

Throughout, the distance between two measurable supports
\(A,B\subset\Omega\) is measured by the Lebesgue measure of their symmetric
difference,
denoted by \(|A\triangle B|\).

Equivalently, we can write  $|A\triangle B|=\norm{\chi_A - \chi_B}_{L^1(\Omega)}$.
Hence \(|A\triangle B|\to0\) means that the two actuator supports differ only on a set of small measure.
This is the natural notion of closeness for the localized perturbations considered below.
\begin{lemma}
\label{lem:eta_chi_identity}
Let \(A,B\subset\Omega\) be measurable sets and define  $\eta_{A,B}:=\chi_B-\chi_A$.
Then it holds
\[
    2\eta_{A,B}\chi_A+\eta_{A,B}^2
    =
    \eta_{A,B} ( \chi_A + \chi_B)
    =
    \eta_{A,B}.
\]
\end{lemma}

\begin{proof}
The claim follows from the identities
\(\chi_A^2=\chi_A\), \(\chi_B^2=\chi_B\), and \(\chi_A\chi_B=\chi_B\chi_A\).
\end{proof}

\subsection{Solvability of closed-loop Lyapunov equation}

The following estimate is the Lyapunov analogue of the standard inverse
generator formula for exponentially stable semigroups; see~\cite[Thm.~II.1.10(i)]{EngelNagel2000}. For the related observability Gramian and Lyapunov equation in semigroup
control theory, see also~\cite[Sec.~5.1, Thm.~5.1.1]{TucsnakWeiss2009}. We include the short proof for completeness.

\begin{lemma}
\label{lem:standard_lyapunov_inverse}
Let $A\in \Aad$ so that
\(\A_A^{\cl}\) generates an exponentially stable \(C_0\)-semi\-group
\(S_A^{\cl}(t)\) on \(H\).
Define
\(\LA X:=(\A_A^{\cl})^*X+X\A_A^{\cl}\).
Then, for every \(F\in\Ls(H)\),
\[
    \LA^{-1}F
    :=
    -\int_0^\infty
    S_A^{\cl}(t)^*F S_A^{\cl}(t)\,\dd t
\]
is well-defined in \(\Ls(H)\) and solves \(\LA X=F\) in the weak Lyapunov
sense.
Moreover, there exists a constant \(C>0\) such that
\begin{equation}
\label{eq:lyapunov_inverse_estimate}
    \norm{\LA^{-1}F}_{\Ls(H)}
    \leq
    C\norm{F}_{\Ls(H)}.
\end{equation}
\end{lemma}

\begin{proof}
By exponential stability, \eqref{eq:closed_loop_exp_stability} holds with
constants \(M_A\geq1\) and \(\gamma_A>0\). For \(F\in\Ls(H)\), define
\[
    X:=
    -\int_0^\infty
    S_A^{\cl}(t)^*F S_A^{\cl}(t)\,\dd t .
\]
The integral converges in \(\Ls(H)\), since
\[
\begin{aligned}
    \int_0^\infty
    \norm{S_A^{\cl}(t)^*F S_A^{\cl}(t)}_{\Ls(H)}
    \,\dd t
    &\leq
    M_A^2\norm{F}_{\Ls(H)}
    \int_0^\infty e^{-2\gamma_A t}\,\dd t
    <\infty .
\end{aligned}
\]
Thus \(X\in\Ls(H)\), and the same estimate gives $\norm{X}_{\Ls(H)} \leq C\norm{F}_{\Ls(H)}$.
This proves \eqref{eq:lyapunov_inverse_estimate}. It remains to verify the Lyapunov identity \(\LA X=F\). Since \(\A_A^{\cl}\) is generally unbounded, the Lyapunov equation is understood in the weak sense on
\(D(\A_A^{\cl})\).
% The integral nevertheless defines a bounded operator on all of $H$.
For \(T>0\) we set
\[
    X_T:=
    -\int_0^T
    S_A^{\cl}(t)^*F S_A^{\cl}(t)\,\dd t ,
\]
so that $X_T \to X$ in $\Ls(H)$ for $T\to\infty$.
Let \(z,w\in D(\A_A^{\cl})\) be arbitrary. Differentiating the scalar function
\(t\mapsto (F S_A^{\cl}(t)z,S_A^{\cl}(t)w)_H\) and integrating over
\((0,T)\), one obtains
\[
\begin{aligned}
    (X_T\A_A^{\cl}z,w)_H
    +
    (X_Tz,\A_A^{\cl}w)_H
    =
    (Fz,w)_H
    -
    (S_A^{\cl}(T)^*F S_A^{\cl}(T)z,w)_H .
\end{aligned}
\]
Together with the fact that
\[
    \norm{S_A^{\cl}(T)^*F S_A^{\cl}(T)}_{\Ls(H)}
    \leq
    M_A^2e^{-2\gamma_A T}\norm{F}_{\Ls(H)}
    \to0
    \qquad
    \text{as }T\to\infty,
\]
we may pass to the limit and obtain
\[
    (X\A_A^{\cl}z,w)_H+(Xz,\A_A^{\cl}w)_H=(Fz,w)_H,
    \qquad z,w\in D(\A_A^{\cl}).
\]
Hence \(X\) solves \(\LA X=F\) in the weak Lyapunov sense, and therefore
\(X=\LA^{-1}F\).
\end{proof}

\subsection{Support perturbation in Riccati equation}

To prove the localized Riccati sensitivity estimate, we first rewrite the
Riccati equation for the support  \(B\) as a perturbation of the Riccati equation for the
reference support \(A\). The following identity isolates the leading
localized term, which depends only on the support variation
\(\eta_{A,B}=\chi_B-\chi_A\), and separates it from the terms containing
\(\Delta\Pi=\Pi_B-\Pi_A\).

\begin{lemma}
\label{lem:riccati_perturbation_identity}
Let \(A,B\in\Aad\) be stabilizing supports and set
\(\Delta\Pi:=\Pi_B-\Pi_A\) and
\(\eta_{A,B}:=\chi_B-\chi_A\). Then \(\Delta\Pi\) satisfies, in the weak
Lyapunov sense,
\begin{equation}
\label{eq:riccati_perturbation_identity}
\begin{aligned}
    &
    (\A_A^{\cl})^*\Delta\Pi
    +
    \Delta\Pi\A_A^{\cl}
    \\
    &=
    \frac1\alpha
    \Big[
        \Pi_A\eta_{A,B}\Pi_A
        +
        \Pi_A\eta_{A,B}\Delta\Pi
        +
        \Delta\Pi\eta_{A,B}\Pi_A
        +
        \Delta\Pi\chi_A\Delta\Pi
        +
        \Delta\Pi\eta_{A,B}\Delta\Pi
    \Big].
\end{aligned}
\end{equation}
\end{lemma}

\begin{proof}
Subtracting the Riccati equation for \(A\) from the one for \(B\) gives
\[
    \A^*\Delta\Pi+\Delta\Pi\A
    =
    \frac1\alpha
    \left(
        \Pi_B\chi_B\Pi_B-\Pi_A\chi_A\Pi_A
    \right).
\]
Using \(\Pi_B=\Pi_A+\Delta\Pi\) and
\(\chi_B=\chi_A+\eta_{A,B}\), we obtain
\[
\begin{aligned}
    \Pi_B\chi_B\Pi_B-\Pi_A\chi_A\Pi_A
    &=
    \Pi_A\chi_A\Delta\Pi
    +
    \Delta\Pi\chi_A\Pi_A
    +
    \Pi_A\eta_{A,B}\Pi_A
    \\
    &\quad+
    \Pi_A\eta_{A,B}\Delta\Pi
    +
    \Delta\Pi\eta_{A,B}\Pi_A
    +
    \Delta\Pi\chi_A\Delta\Pi
    +
    \Delta\Pi\eta_{A,B}\Delta\Pi .
\end{aligned}
\]
Moreover, since
\(\A_A^{\cl}=\A-\alpha^{-1}\chi_A\Pi_A\) and
\((\A_A^{\cl})^*=\A^*-\alpha^{-1}\Pi_A\chi_A\), we have
\[
\begin{aligned}
    &
    (\A_A^{\cl})^*\Delta\Pi
    +
    \Delta\Pi\A_A^{\cl}
    =
    \A^*\Delta\Pi+\Delta\Pi\A
    -
    \frac1\alpha\Pi_A\chi_A\Delta\Pi
    -
    \frac1\alpha\Delta\Pi\chi_A\Pi_A .
\end{aligned}
\]
Thus the two terms that are linear in \(\Delta\Pi\) and contain \(\chi_A\) are
absorbed into the closed-loop left-hand side, which gives
\eqref{eq:riccati_perturbation_identity}.
\end{proof}

\section{Riccati regularity and sensitivity with respect to the support} \label{sec:riccati_sensitivity}

This section contains the main analytical estimates behind the support sensitivity result.

\subsection{Riccati smoothing and feedback regularity} \label{subsec:riccati_smoothing}

The following regularity result is used to justify weak pairings, to prove
local stabilizability, and to motivate the stronger regularity hypotheses used
later in the localized sensitivity analysis.

\begin{proposition}[Riccati smoothing]
\label{prop:riccati_smoothing}
Let \(A\in\Aad\)
so that
the closed-loop generator \(\A_A^{\cl}\)
generates an exponentially stable analytic semigroup \(S_A^{\cl}(t)\) on \(H\).
Then, for every
\(0\leq\theta<1\), we have
\[
    \Pi_A\in
    \Ls\left(
        H,
        D\left((-\A_A^{\cl,*})^\theta\right)
    \right).
\]
% \red
% In particular,
% \(\Pi_A\in\Ls(H,V)\)
% ??
% \der%
\end{proposition}

\begin{proof}
% We first explain why the closed-loop semigroup  is analytic.
The operator $\A$ generates an analytic semigroup.
Since \(\chi_A\KA\in\Ls(H)\), the bounded
perturbation theorem for analytic semigroups \cite[Ch.~3, Cor.~2.2, p.~81]{Pazy1983}
implies that
the closed-loop generator \(\A_A^{\cl}=\A+\chi_A\KA\) also generates an analytic semigroup \(S_A^{\cl}(t)\).

We set \(L_A:=Q^*Q+\alpha\KA^*\KA\). Then \(L_A\in\Ls(H)\), and by the
closed-loop Lyapunov form \eqref{eq:riccati_operator2} together
with \cref{lem:standard_lyapunov_inverse}, the stabilizing Riccati
operator admits the closed-loop representation
\begin{equation}
\label{eq:riccati_closed_loop_rep}
    \Pi_A
    =
    \int_0^\infty
    S_A^{\cl}(t)^*
    L_A
    S_A^{\cl}(t)
    \,\dd t.
\end{equation}
Let \(0\leq\theta<1\). Since \(S_A^{\cl}(t)^*\) is analytic and exponentially
stable, the standard fractional-power estimate gives, possibly after decreasing \(\gamma_A\) in \eqref{eq:closed_loop_exp_stability},
\begin{equation}
\label{eq:fractional_semigroup_estimate}
    \norm{
        (-\A_A^{\cl,*})^\theta
        S_A^{\cl}(t)^*
    }_{\Ls(H)}
    \leq
    C_\theta t^{-\theta}e^{-\gamma_A t},
    \qquad
    t>0.
\end{equation}
See ~\cite[Ch.~2, Thm.~6.13, p.~74]{Pazy1983}. Moreover, using
\eqref{eq:closed_loop_exp_stability}, we obtain

\begin{equation*}
    \begin{aligned}
        & \norm{
        (-\A_A^{\cl,*})^\theta
        S_A^{\cl}(t)^*
        L_A
        S_A^{\cl}(t)
    }_{\Ls(H)}
    \\
        & \qquad \leq
    \norm{
        (-\A_A^{\cl,*})^\theta
        S_A^{\cl}(t)^*
    }_{\Ls(H)}
    \norm{L_A}_{\Ls(H)}
    \norm{S_A^{\cl}(t)}_{\Ls(H)}
    \leq
    C t^{-\theta}e^{-2\gamma_A t}.
    \end{aligned}
\end{equation*}
% Since \(t\mapsto t^{-\theta}e^{-2\gamma_A t}\) is integrable on
% \((0,\infty)\) for \(\theta<1\)
Hence, the map
\(t\mapsto (-\A_A^{\cl,*})^\theta S_A^{\cl}(t)^*L_AS_A^{\cl}(t)\) is Bochner integrable in \(\Ls(H)\).

We now justify applying the fractional power to the integral in
\eqref{eq:riccati_closed_loop_rep}. For every \(t>0\), analyticity gives
\(S_A^{\cl}(t)^*H\subset D((-\A_A^{\cl,*})^\theta)\), and the operator
\((-\A_A^{\cl,*})^\theta\) is closed; see
Pazy~\cite[Ch.~2, Sec.~6]{Pazy1983}. Since for $y\in H$ both
\(t\mapsto S_A^{\cl}(t)^*L_AS_A^{\cl}(t)y\) and
\(t\mapsto (-\A_A^{\cl,*})^\theta S_A^{\cl}(t)^*L_AS_A^{\cl}(t)y\) are
Bochner integrable, closedness implies
\[
    \Pi_Ay
    =
    \int_0^\infty
    S_A^{\cl}(t)^*L_AS_A^{\cl}(t)y
    \,\dd t
    \in
    D((-\A_A^{\cl,*})^\theta)
\]
and
\[
    (-\A_A^{\cl,*})^\theta\Pi_Ay
    =
    \int_0^\infty
    (-\A_A^{\cl,*})^\theta
    S_A^{\cl}(t)^*L_AS_A^{\cl}(t)y
    \,\dd t .
\]
Therefore, we have
\[
\begin{aligned}
    \norm{(-\A_A^{\cl,*})^\theta\Pi_Ay}_H
    &\leq
    \int_0^\infty
    \norm{
        (-\A_A^{\cl,*})^\theta
        S_A^{\cl}(t)^*
        L_A
        S_A^{\cl}(t)y
    }_H
    \,\dd t
    \\
    &\leq
    C\norm{y}_H
    \int_0^\infty
    t^{-\theta}e^{-2\gamma_A t}
    \,\dd t
    \leq
    C_\theta\norm{y}_H ,
\end{aligned}
\]
which  yields $\norm{\Pi_Ay}_{D((-\A_A^{\cl,*})^\theta)} \leq C_\theta\norm{y}_H$.
This proves the asserted fractional-domain smoothing.
\end{proof}

\begin{corollary}
\label{cor:spatial_consequences_riccati_smoothing}
Let the assumptions of \cref{prop:riccati_smoothing} hold.
Let $\A$ be a second-order differential operator as in \cref{prop:parabolic_is_nice}.
Then
the following consequences are valid.
\begin{enumerate}[label=\textup{(C\arabic*)}]
    \item \label{cor:C1_fractional_embeddings}
    It holds
    \(\Pi_A\in\Ls(H,V)\). Moreover, in spatial dimensions \(d\leq3\), one also
    obtains \(\Pi_A\in\Ls(H,L^\infty(\Omega))\).

    \item \label{cor:C2_Pi_Lq_to_H}
    If \(\Pi_A\in\Ls(H,L^\infty(\Omega))\), then \(\Pi_A\) extends uniquely
    to a bounded operator \(\Pi_A:L^1(\Omega)\to H\). Consequently,
    \(\Pi_A:L^q(\Omega)\to H\) for every \(1\leq q\leq2\).

    \item \label{cor:C3_feedback_Lp}
    If \(\Pi_A\in\Ls(H,V)\) and \(V\hookrightarrow L^p(\Omega)\) for some
    \(p>2\), then \(\KA\in\Ls(H,L^p(\Omega))\).
\end{enumerate}

\end{corollary}

\begin{proof}
We first prove \ref{cor:C1_fractional_embeddings}. \Cref{prop:riccati_smoothing} gives
\(\Pi_A\in\Ls(H,D((-\A_A^{\cl,*})^\theta))\) for every \(0\leq\theta<1\). For the underlying second-order elliptic operator, the spatial regularity
can be obtained directly from its form. Since the closed-loop term is
bounded on \(H\), the adjoint closed-loop form retains the G{\aa}rding
inequality. Together with analyticity, this gives
\(\norm{S_A^{\cl}(t)^*}_{\Ls(H,V)}\leq C_A t^{-1/2}\) for \(0<t\leq1\).
Together with exponential stability, this estimate and
\eqref{eq:riccati_closed_loop_rep} give \(\Pi_A\in\Ls(H,V)\).
Similarly, the smoothing estimate used in \cref{lem:Riccati_linfty_estimate}
has order \(t^{-d/4}\), which is integrable at zero for \(d\leq3\), and gives
\(\Pi_A\in\Ls(H,L^\infty(\Omega))\). Thus both inclusions in
\ref{cor:C1_fractional_embeddings} follow.

We prove \ref{cor:C2_Pi_Lq_to_H}. %Let \(f\in L^1(\Omega)\cap H\).
Since \(\Omega\) is bounded, \(H=L^2(\Omega)\) is dense in \(L^1(\Omega)\).
%; thus this is a dense subspace on which \(\Pi_A f\) is already defined.
Using \(\Pi_A\in\Ls(H, L^\infty(\Omega))\) and self-adjointness of $\Pi_A$ on $H$, we obtain for $f\in H$
\[
\norm{\Pi_Af}_H = \sup_{ \norm{h}_H\le 1} \ip{\Pi_A f}{h}_H
    \leq
    \norm{f}_{L^1(\Omega)}
    \norm{\Pi_A}_{\Ls(H,L^\infty(\Omega))}
    \leq
    C\norm{f}_{L^1(\Omega)}.
\]
Since \(H\) is dense in
\(L^1(\Omega)\), the operator \(\Pi_A\) can be uniquely extended
to a bounded
operator from \(L^1(\Omega)\) to \(H\). By interpolation between \(L^1(\Omega)\) and \(L^2(\Omega)\), see~\cite[Sec.~1.2]{Lunardi1995}, we obtain
\(\Pi_A:L^q(\Omega)\to H\) for every \(1\leq q\leq2\).
Finally, \ref{cor:C3_feedback_Lp} follows from \(\KA=-\alpha^{-1}\chi_A\Pi_A\). %Indeed,
% if \(V\hookrightarrow L^p(\Omega)\), then
% \[
%     \norm{\KA y}_{L^p(\Omega)}
%     \leq
%     \frac1\alpha\norm{\Pi_Ay}_{L^p(\Omega)}
%     \leq
%     C\norm{\Pi_Ay}_V
%     \leq
%     C\norm{y}_H .
% \]
\end{proof}

\subsection{Local preservation of stabilizability} \label{subsec:local_stabilizability}

The following lemma shows that the set of stabilizing supports is locally open
with respect to the symmetric-difference metric introduced in
\cref{subsec:support_perturbations}. This is needed because a small topological
perturbation of a stabilizing support should still define a well-posed
infinite-horizon LQR problem.
% The proof uses the elementary localization
% estimates from \cref{lem:localization_symmetric_difference,lem:localized_Lp_L2}.
The assumption \(\KA\in\Ls(H,L^p(\Omega))\) for some $p>2$ is justified by \cref{cor:spatial_consequences_riccati_smoothing}.

\begin{lemma}[Local preservation of stabilizability]
\label{lem:local_stabilizability}
Let \(A\in\Aad\).
Suppose that the feedback operator satisfies \(\KA\in\Ls(H,L^p(\Omega))\) for some $p>2$.
Then there exists
\(\delta_A>0\) such that for every measurable \(B\subset\Omega\) with
\(|A\triangle B|<\delta_A\), the pair \((\A,\chi_B)\) is exponentially
stabilizable.

More precisely, the feedback \(\KA\) exponentially stabilizes
the system with support \(B\). Denote by \(S_B^A(t)\) the semigroup generated
by \(\A+\chi_B\KA\), then
\[
    \norm{S_B^A(t)}_{\Ls(H)}
    \leq
    M_A
    \exp\left(
        -\left[
            \gamma_A
            -
            M_A C_A |A\triangle B|^{\frac12-\frac1p}
        \right]t
    \right),
    \qquad
    t\geq0,
\]
where \(C_A:=\norm{\KA}_{\Ls(H,L^p(\Omega))}\).
Thus, choosing \(\delta_A>0\) sufficiently small so that
\(M_A C_A |A\triangle B|^{\frac12-\frac1p}<\gamma_A\), the semigroup
\(S_B^A(t)\) is exponentially stable.
\end{lemma}

\begin{proof}
We use the feedback \(\KA\), which stabilizes the reference support \(A\), as
a test feedback for the perturbed support \(B\). Setting $E_{A,B}:=(\chi_B-\chi_A)\KA$, we can  write  $\A+\chi_B\KA =\A_A^{\cl}+E_{A,B}$. We first estimate the perturbation \(E_{A,B}\).
% Using \eqref{eq:local_pertub2} from \cref{lem:localization_symmetric_difference} and using \cref{lem:localized_Lp_L2},
For   \(y\in H\) we get
\begin{equation}
\label{eq:lem_local0}
\norm{E_{A,B}y}_H =\norm{(\chi_B-\chi_A)\KA y}_{L^2(\Omega)} = \norm{\KA y}_{L^2(A\triangle B)} \leq  |A\triangle B|^{\frac12-\frac1p} \norm{\KA y}_{L^p(\Omega)}
\end{equation}
Due to \(\KA\in\Ls(H,L^p(\Omega))\), it follows 
\begin{equation}
\label{eq:EAB_bound}
    \norm{E_{A,B}}_{\Ls(H)}
    \leq
    C_A|A\triangle B|^{\frac12-\frac1p}
\end{equation}
with $C_A:=\norm{\KA}_{\Ls(H,L^p(\Omega))}$.
Since \(p>2\), the exponent $\frac12-\frac1p$ is positive, and therefore $\norm{E_{A,B}}_{\Ls(H)}\to0$ as  $|A\triangle B|\to0$.

By the bounded perturbation theorem for \(C_0\)-semigroups, see \cite[Ch.~3, Thm.~1.1, p.~76]{Pazy1983},
the operator $\A_A^{\cl}+E_{A,B}= \A+\chi_B\KA$ generates a \(C_0\)-semigroup on \(H\), which we denote by \(S_B^A(t)\). We now show that this semigroup remains exponentially stable if \(|A\triangle B|\) is small enough. By the variation-of-constants formula, we can write
\[
    S_B^A(t)
    =
    S_A^{\cl}(t)
    +
    \int_0^t
    S_A^{\cl}(t-s)E_{A,B}S_B^A(s)
    \,\dd s .
\]
Taking operator norms and using the exponential stability of \(S_A^{\cl}(t)\), we
obtain
\[
\begin{aligned}
    \norm{S_B^A(t)}_{\Ls(H)}
    &\leq
    M_Ae^{-\gamma_A t}
    +
    \int_0^t
    M_Ae^{-\gamma_A(t-s)}
    \norm{E_{A,B}}_{\Ls(H)}
    \norm{S_B^A(s)}_{\Ls(H)}
    \,\dd s .
\end{aligned}
\]
Multiplying by \(e^{\gamma_A t}\) and applying Gronwall's inequality gives
\[
    e^{\gamma_A t}\norm{S_B^A(t)}_{\Ls(H)}
    \leq
    M_A
    \exp\left(
        M_A\norm{E_{A,B}}_{\Ls(H)}t
    \right).
\]
Therefore, using \eqref{eq:EAB_bound}, we obtain for every $t\geq0$ that
\begin{equation}
\begin{aligned}
 \norm{S_B^A(t)}_{\Ls(H)} &\leq M_A \exp\left(-\left[\gamma_A-M_A\norm{E_{A,B}}_{\Ls(H)} \right]t \right)\\ &\leq M_A \exp\left(-\left[\gamma_A-M_A C_A |A\triangle B|^{\frac12-\frac1p}\right]t \right).
\end{aligned}
\end{equation}
Let us choose \(\delta_A>0\) such that  $M_A C_A\delta_A^{\frac12-\frac1p} <\gamma_A$. Then, whenever $|A\triangle B|<\delta_A$ holds, the decay rate
 $\gamma_A  -M_A C_A |A\triangle B|^{\frac12-\frac1p}$ is strictly positive. Hence \(\A+\chi_B\KA\) generates an exponentially stable semigroup. Thus the feedback \(\KA\) stabilizes the system with support \(B\). Therefore
the pair \((\A,\chi_B)\) is exponentially stabilizable.
\end{proof}

\subsection{Comparison of value functions} \label{subsec:value_comparison}

We now derive the value function comparison identity that will be used to obtain the
first-order expansion. This identity is the key step where the difference of
two support-dependent value functions is written along a single reference
trajectory.

\begin{proposition}
\label{prop:value_comparison_identity}
Let \(A,B\in\Aad\) be admissible stabilizing supports. Let \(y_A\) be the
optimal closed-loop state corresponding to \(A\), and let
\(\Pi_A,\Pi_B\) and \(\KA,\KB\) denote the corresponding Riccati operators and
feedback operators. Assume that the smoothing property
\(\Pi_B\in\Ls(H,V)\) holds.
Then it holds
\begin{equation}
\label{eq:value_difference_identity}
\begin{aligned}
    V_A(y_0)-V_B(y_0)
    &=
    2\int_0^\infty
    \left(
        (\chi_A-\chi_B)\KA y_A(t),
        \Pi_B y_A(t)
    \right)_H
    \dd t
    \\
    &\quad+
    \alpha
    \norm{(\KB-\KA)y_A}_{L^2(0,\infty;H)}^2 .
\end{aligned}
\end{equation}
\end{proposition}

\begin{proof}
By \cref{lem:closed_loop_W_regular}, the reference closed-loop state
satisfies $y_A\in W(0,\infty)$. Moreover, since \(\Pi_B\in\Ls(H,V)\) and
\(y_A\in L^2(0,\infty;H)\), we have $\Pi_B y_A\in L^2(0,\infty;V)$. Hence the duality pairing $ \duality{\dot y_A(t)}{\Pi_B y_A(t)}_{V',V}$ is integrable on \((0,\infty)\).

Since \(y_A\) is the optimal state associated with the support \(A\), we have
\[
    V_A(y_0)
    =
    \int_0^\infty
    \left(
        \norm{Qy_A(t)}_{H_Q}^2
        +
        \alpha\norm{\KA y_A(t)}_H^2
    \right)
    \dd t.
\]
The Riccati identity associated with the support \(B\) gives us for sufficiently
regular \(z\) that
 \[
    2\duality{\A z}{\Pi_B z}_{V',V}
    -
    \alpha\norm{\KB z}_H^2
    +
    \norm{Qz}_{H_Q}^2
    =
    0.
\]
Using this identity with \(z=y_A(t)\), and using the closed-loop equation
\[
    \dot y_A(t)=\A y_A(t)+\chi_A\KA y_A(t),
\]
we can write that
\begin{multline*}
    \norm{Qy_A(t)}_{H_Q}^2
     =
    -2\duality{\A y_A(t)}{\Pi_B y_A(t)}_{V',V}
    +
    \alpha\norm{\KB y_A(t)}_H^2
    \\
    =
    -2\duality{\dot y_A(t)}{\Pi_B y_A(t)}_{V',V}
    +
    2\left(
        \chi_A\KA y_A(t),
        \Pi_B y_A(t)
    \right)_H+ \alpha\norm{\KB y_A(t)}_H^2.
\end{multline*}
Therefore, we have
\begin{multline}
\label{eq:lem12_1}
    \norm{Qy_A(t)}_{H_Q}^2
    +
    \alpha\norm{\KA y_A(t)}_H^2
    =
    -2\duality{\dot y_A(t)}{\Pi_B y_A(t)}_{V',V}
    \\+
    2\left(
        \chi_A\KA y_A(t),
        \Pi_B y_A(t)
    \right)_H+
    \alpha\norm{\KB y_A(t)}_H^2
    +
    \alpha\norm{\KA y_A(t)}_H^2.
\end{multline}
Using $\KB=-\frac1\alpha\chi_B\Pi_B$, we have
\[
    \alpha
    \left(
        \KA y_A(t),
        \KB y_A(t)
    \right)_H
    =
    -
    \left(
        \chi_B\KA y_A(t),
        \Pi_B y_A(t)
    \right)_H.
\]
Consequently, we have
\begin{multline*}
    \alpha\norm{\KB y_A(t)}_H^2
    +
    \alpha\norm{\KA y_A(t)}_H^2
    +
    2\left(
        \chi_A\KA y_A(t),
        \Pi_B y_A(t)
    \right)_H
    \\
    =
    \alpha\norm{(\KB-\KA)y_A(t)}_H^2
    +
    2\left(
        (\chi_A-\chi_B)\KA y_A(t),
        \Pi_B y_A(t)
    \right)_H.
\end{multline*}
This together with \eqref{eq:lem12_1} gives us
\begin{multline*}
    \norm{Qy_A(t)}_{H_Q}^2
    +
    \alpha\norm{\KA y_A(t)}_H^2
    =
    -2\duality{\dot y_A(t)}{\Pi_B y_A(t)}_{V',V}
    \\
    +
    2\left(
        (\chi_A-\chi_B)\KA y_A(t),
        \Pi_B y_A(t)
    \right)_H+
    \alpha\norm{(\KB-\KA)y_A(t)}_H^2.
\end{multline*}
Integrating over \((0,T)\), we get that
\begin{equation}
\label{eq:lem12_2}
\begin{aligned}
  &
    \int_0^T
    \left(
        \norm{Qy_A(t)}_{H_Q}^2
        +
        \alpha\norm{\KA y_A(t)}_H^2
    \right)
    \dd t=
    -2\int_0^T
    \duality{\dot y_A(t)}{\Pi_B y_A(t)}_{V',V}
    \dd t
    \\
    &\quad+
    2\int_0^T
    \left(
        (\chi_A-\chi_B)\KA y_A(t),
        \Pi_B y_A(t)
    \right)_H
    \dd t +
    \alpha
    \norm{(\KB-\KA)y_A}_{L^2(0,T;H)}^2.
\end{aligned}
\end{equation}
By the chain rule, we have for almost every \(t>0\) that
\[
    \frac{\dd}{\dd t}
    \ip{\Pi_B y_A(t)}{y_A(t)}_H
    =
    2\duality{\dot y_A(t)}{\Pi_B y_A(t)}_{V',V}.
\]
Hence, we can write
\[
\begin{aligned}
    -2\int_0^T
    \duality{\dot y_A(t)}{\Pi_B y_A(t)}_{V',V}
    \dd t
    &=
    \ip{\Pi_B y_0}{y_0}_H
    -
    \ip{\Pi_B y_A(T)}{y_A(T)}_H.
\end{aligned}
\]
Using the fact that $\ip{\Pi_B y_A(T)}{y_A(T)}_H\to0$ as  $T\to\infty$ and letting \(T\to\infty\) in \eqref{eq:lem12_2}, and using $V_B(y_0)=\ip{\Pi_B y_0}{y_0}_H$, we obtain
\[
\begin{aligned}
    V_A(y_0)
    &=
    V_B(y_0)
    +
    2\int_0^\infty
    \left(
        (\chi_A-\chi_B)\KA y_A(t),
        \Pi_B y_A(t)
    \right)_H
    \dd t
    \\
    &\quad+
    \alpha
    \norm{(\KB-\KA)y_A}_{L^2(0,\infty;H)}^2.
\end{aligned}
\]
This is exactly \eqref{eq:value_difference_identity}.
\end{proof}

Using similar arguments, we find the following one-sided estimate, which is analogous to \cite[Lem.~3.5]{Wachsmuth2024}.

\begin{proposition}
\label{prop:value_one_sided_estimate}
Let \(A,B\in\Aad\) be admissible stabilizing supports. Let \(y_A\) be the
optimal closed-loop state corresponding to \(A\), and let
\(\Pi_A,\Pi_B\) denote the corresponding Riccati operators.
Assume that the smoothing property
\(\Pi_A,\Pi_B\in\Ls(H,V)\) holds. Then it holds
\begin{equation}
\label{eq:value_one_sided_estimate}
  V_B(y_0)-V_A(y_0)
  \ge -\frac1\alpha
    \int_\Omega
    (\chi_B-\chi_A)(x)
    \int_0^\infty |\Pi_Ay_A(t)(x)|^2\dd t
    \dd x.
\end{equation}
\end{proposition}
\begin{proof}
 Let $y_A, y_B$  and $u_A,u_B$ be the optimal closed loop state and controls to $A,B$, respectively.
 Set $\delta y := y_B-y_A$ and $\delta u := u_B-u_A$.
 Due to convexity, we have
 \[
  V_B(y_0)-V_A(y_0) \ge 2\int_0^\infty \ip{ Qy_A(t) }{ Q\delta y(t)}_{H_Q} + \alpha \ip{ u_A(t)}{\delta u(t)}_H + \frac\alpha2 \norm{ \delta u(t)}_H^2 \dd t.
 \]
Using the Riccati equation for $\Pi_A$ together with $\dot  y_A = \A y_A -\frac1\alpha\chi_A\Pi_Ay_A$, yields
\[
 0 = \ip{ Qy_A (t)}{ Q\delta y(t)}_{H_Q} + \duality{\dot y_A(t)  }{\Pi_A \delta y(t)}_{V',V}
 + \duality{\A \delta y(t)  }{\Pi_A y_A(t)}_{V',V}.
%  0 = \ip{ Qy_A (t)}{ Q\delta y(t)}_{H_Q} + \duality{\A y_A(t)  }{\Pi_A \delta y(t)}_{V',V}
%  + \duality{\A \delta y(t)  }{\Pi_A y_A(t)}_{V',V} + \ip{ \KA y_A(t)}{\Pi_A \delta y(t)}_H.
\]
Using the closed loop equations for $y_A, y_B$, we find $\A \delta y = \delta \dot y - \delta u$, so that we obtain
\begin{multline*}
 0=\ip{ Qy_A(t) }{ Q\delta y(t)}_{H_Q} + \duality{\dot y_A(t)}{\Pi_A \delta y(t)}_{V',V} \\
 + \duality{\delta \dot y(t)  }{\Pi_A y_A(t)}_{V',V}
 -\ip{ \delta u(t)}{\Pi_A y_A(t)}_H.
\end{multline*}
Since
$
 \duality{\dot y_A(t)}{\Pi_A \delta y(t)}_{V',V}
 + \duality{\delta \dot y(t)  }{\Pi_A y_A(t)}_{V',V} = \frac\dd{\dd t} \ip{ y_A(t) }{ \Pi_A \delta y(t)}_H,
$
we find
\begin{multline*}
 \int_0^\infty  \duality{\dot y_A(t)}{\Pi_A \delta y(t)}_{V',V}
 + \duality{\delta \dot y (t) }{\Pi_A y_A(t)}_{V',V} \dd t \\
 = \lim_{T\to\infty} \ip{ y_A(T) }{ \Pi_A \delta y(T)}_H-\ip{ y_0 }{ \Pi_A \delta y(0)}_H=0.
\end{multline*}
This implies
\[
 \int_0^\infty \ip{ Qy_A(t) }{ Q\delta y(t)}_{H_Q} \dd t = \int_0^\infty \ip{ \delta u(t)}{\Pi_A y_A(t)}_H \dd t,
\]
which results in the inequality
 \[\begin{aligned}
  V_B(y_0)-V_A(y_0) &\ge 2\int_0^\infty \ip{ \delta u(t)}{\Pi_A y_A(t)}_H + \alpha \ip{ u_A(t)}{\delta u(t)}_H
  +\frac\alpha2 \norm{ \delta u(t)}_H^2 \,\dd t\\
  &=2\int_0^\infty \ip{ \delta u(t)}{ (1 - \chi_A) \Pi_A y_A(t)}_H + \frac\alpha2 \norm{ \delta u(t)}_H^2 \dd t.
 \end{aligned}\]
We split the integral on $\Omega$, which corresponds to the inner product on $H$, into
contributions on $A\cap B$, $A\setminus B$, $B\setminus A$. Note that on $\Omega \setminus (A \cup B)$ we have $\delta u=0$.
For the contribution on $A\cap B$ we find
\[
 \int_0^\infty \int_{A\cap B}  \delta u(t,x) (1 - \chi_A(x)) \Pi_A y_A(t)(x) + \frac\alpha2 |\delta u(t,x)|^2 \dd x \dd t
 %= \int_0^\infty \int_{A\cap B}\frac\alpha2 |\delta u(t,x)|^2 \dd x \dd t
 \ge 0.
\]
Using $\min_{v\in \R} \left( v\cdot p + \frac\alpha2 v^2 \right)= -\frac1{2\alpha}p^2$, we get on $B\setminus A$
\begin{multline*}
 \int_0^\infty \int_{B\setminus A}  \delta u(t,x) (1 - \chi_A(x)) \Pi_A y_A(t)(x) + \frac\alpha2 |\delta u(t,x)|^2 \dd x \dd t\\
% = \int_0^\infty \int_{B\setminus A}  \delta u(t,x) \Pi_A y_A(t)(x) + \frac\alpha2 |\delta u(t,x)|^2 \dd x \dd t
 \ge \int_0^\infty \int_{B\setminus A} -\frac1{2\alpha} |\Pi_A y_A(t)(x) |^2  \dd x \dd t.
\end{multline*}
For $x\in A\setminus B$ it holds $\delta u(t,x) = -u_A = \frac1\alpha \Pi_A y_A(t)(x)$, so that
\begin{multline*}
 \int_0^\infty \int_{A\setminus B}  \delta u(t,x) (1 - \chi_A(x)) \Pi_A y_A(t)(x) + \frac\alpha2 |\delta u(t,x)|^2 \dd x \dd t\\
= \int_0^\infty \int_{A\setminus B}   \frac\alpha2 |\delta u(t,x)|^2 \dd x \dd t
 = \int_0^\infty \int_{A\setminus B}  \frac1{2\alpha} |\Pi_A y_A(t)(x) |^2  \dd x \dd t.
\end{multline*}
Collecting these estimates, we get
\[
\begin{aligned}
  V_B(y_0)-V_A(y_0) &\ge 2\int_0^\infty \ip{ \delta u}{ (1 - \chi_A) \Pi_A y_A(t)}_H + \frac\alpha2 \norm{ \delta u}_H^2 \dd t\\
  & \ge -\frac1\alpha \int_0^\infty \int_\Omega (\chi_B-\chi_A)(x)  |\Pi_A y_A(t)(x) |^2 \dd x  \dd t,
 \end{aligned}
\]
which proves the claim.
\end{proof}

\subsection{Exact expansion and remainder estimate} \label{subsec:exact_expansion}
We now use \cref{prop:value_comparison_identity} to derive an exact
expansion of \(J(B)-J(A)\), with an explicit remainder. This is the point where the comparison identity is connected to the topological derivative.
For two stabilizing supports \(A\) and \(B\), we use the notation
\[
    \eta_{A,B}:=\chi_B-\chi_A,
    \qquad
    \psi_A(t):=\Pi_Ay_A(t),
    \qquad
    \psi_B(t):=\Pi_By_A(t).
\]
% The elementary identity needed in the extraction of the leading term has
% already been recorded in \cref{lem:eta_chi_identity}.

\begin{lemma}[Exact expansion with explicit remainder]
\label{lem:exact_expansion_remainder}
Let \(A\) and \(B\) be stabilizing supports for which
\cref{prop:value_comparison_identity} applies. Defining
\begin{equation}
\label{eq:switching_function}
    G_A(x)
    :=
    \beta(x)-\frac1\alpha\int_0^\infty |\Pi_Ay_A(t)(x)|^2\dd t.
\end{equation}
Then, we can write
\begin{equation}
\label{eq:exact_expansion_with_remainder}
    J(B)-J(A)
    =
    \int_\Omega \eta_{A,B}(x)G_A(x)\dd x+
    \mathcal R_A(B),
\end{equation}
where
\begin{equation}
\label{eq:remainder_explicit}
% \begin{aligned}
    \mathcal R_A(B)
    :=
    -\frac1\alpha
    \int_0^\infty
    2
    \left(
        \eta_{A,B}\psi_A(t),
        \psi_B(t)-\psi_A(t)
    \right)_H
%     \dd t
%     -
    +
%     \frac1\alpha
%     \int_0^\infty
    \norm{\chi_B(\psi_B(t)-\psi_A(t))}_H^2
    \dd t.
% \end{aligned}
\end{equation}
\end{lemma}

\begin{proof}
Using the fact that $\KA y_A=-\frac1\alpha\chi_A\psi_A$ and   $\KB y_A=-\frac1\alpha\chi_B\psi_B$, we obtain
\begin{equation}
\label{leq:lem13_0}
   \KB y_A-\KA y_A
    =
    -\frac1\alpha\left(\eta_{A,B}\psi_A+
    \chi_B(\psi_B-
    \psi_A)\right).
\end{equation}
 Substituting \eqref{leq:lem13_0} in \eqref{eq:value_difference_identity}, we obtain
\[
\begin{aligned}
    V_B(y_0)-V_A(y_0)
    &=
    -\frac2\alpha
    \int_0^\infty
    \left(
        \eta_{A,B}\chi_A\psi_A(t),
        \psi_B(t)
    \right)_H
    \dd t
    \\
    &\quad
    -
    \frac1\alpha
    \int_0^\infty
    \left\|
        \eta_{A,B}\psi_A(t)
        +
        \chi_B(\psi_B(t)-\psi_A(t))
    \right\|_H^2
    \dd t .
\end{aligned}
\]
We now write $\psi_B(t)=\psi_A(t)+(\psi_B(t)-\psi_A(t))$. Hence the first integral becomes
\[
\begin{aligned}
    &-\frac2\alpha
    \int_0^\infty
    \left(
        \eta_{A,B}\chi_A\psi_A(t),
        \psi_B(t)
    \right)_H
    \dd t
    \\
    &=
    -\frac2\alpha
    \int_0^\infty
    \left(
        \eta_{A,B}\chi_A\psi_A(t),
        \psi_A(t)
    \right)_H
    \dd t  -\frac2\alpha
    \int_0^\infty
    \left(
        \eta_{A,B}\chi_A\psi_A(t),
        \psi_B(t)-\psi_A(t)
    \right)_H
    \dd t .
\end{aligned}
\]
Moreover, expanding the second integral gives
\[
\begin{aligned}
    &
    -\frac1\alpha
    \int_0^\infty
    \left\|
        \eta_{A,B}\psi_A(t)
        +
        \chi_B(\psi_B(t)-\psi_A(t))
    \right\|_H^2
    \dd t =
    -\frac1\alpha
    \int_0^\infty
    \norm{\eta_{A,B}\psi_A(t)}_H^2
    \dd t
    \\
    &\quad
    -\frac2\alpha
    \int_0^\infty
    \left(
        \eta_{A,B}\psi_A(t),
        \chi_B(\psi_B(t)-\psi_A(t))
    \right)_H
    \dd t  -\frac1\alpha
    \int_0^\infty
    \norm{\chi_B(\psi_B(t)-\psi_A(t))}_H^2
    \dd t .
\end{aligned}
\]
The terms containing only \(\psi_A\) are therefore
\[
    -\frac2\alpha
    \int_0^\infty
    \left(
        \eta_{A,B}\chi_A\psi_A(t),
        \psi_A(t)
    \right)_H
    \dd t
    -
    \frac1\alpha
    \int_0^\infty
    \norm{\eta_{A,B}\psi_A(t)}_H^2
    \dd t .
\]
Writing these two terms pointwise, we get
\[
\begin{aligned}
    &-\frac2\alpha
    \int_0^\infty
    \left(
        \eta_{A,B}\chi_A\psi_A(t),
        \psi_A(t)
    \right)_H
    \dd t
    -
    \frac1\alpha
    \int_0^\infty
    \norm{\eta_{A,B}\psi_A(t)}_H^2
    \dd t
    \\
    &=
    -\frac1\alpha
    \int_0^\infty
    \int_\Omega
    \left(
        2\eta_{A,B}(x)\chi_A(x)
        +
        \eta_{A,B}(x)^2
    \right)
    |\psi_A(t,x)|^2
    \dd x\dd t .
\end{aligned}
\]
By \cref{lem:eta_chi_identity},  $2\eta_{A,B}\chi_A+\eta_{A,B}^2=\eta_{A,B}$, and, thus, the pure \(\psi_A\)-terms reduce to
\[
    -\frac1\alpha
    \int_0^\infty
    \int_\Omega
    \eta_{A,B}(x)|\psi_A(t,x)|^2
    \dd x\dd t .
\]
The remaining terms all contain \(\psi_B-\psi_A\), and these terms are
collected in the remainder \(\mathcal R_A(B)\),
where we used $\eta_{A,B}( \chi_A + \chi_B) =  \eta_{A,B}$ from \cref{lem:eta_chi_identity} to simplify
the expression.
Therefore,
\[
    V_B(y_0)-V_A(y_0)
    =
    -\frac1\alpha
    \int_\Omega
    \eta_{A,B}(x)
    \left[
        \int_0^\infty |\psi_A(t,x)|^2\dd t
    \right]
    \dd x+
    \mathcal R_A(B).
\]
Adding  $\int_\Omega(\chi_B-\chi_A)\beta\dd x$ gives \eqref{eq:exact_expansion_with_remainder}.
\end{proof}

\subsection{Localized Riccati sensitivity} \label{subsec:localized_sensitivity}

We now prove a concrete sufficient condition under which the remainder in
\eqref{eq:remainder_explicit} is of higher order. The essential point is that
the localized perturbation  $\eta:=\chi_B-\chi_A$ is supported only on $D:=A\triangle B$.
Thus smallness of \(|A\triangle B|\) can be exploited if the quantities
multiplied by \(\eta\) are sufficiently regular in space.

\begin{theorem}[Localized Riccati sensitivity]
\label{thm:local_branch_and_sensitivity}
Let \(A\in\Aad\) be a stabilizing support
such that \(\Pi_A\in\Ls(H,L^\infty(\Omega))\).
Then there
exist constants \(C_A>0\) and \(\delta_A>0\) such that for every stabilizing
support \(B\in\Aad\) with $|A\triangle B|<\delta_A$, the corresponding Riccati operators satisfy
\begin{equation}
\label{eq:Pi_difference_rate}
    \norm{\Pi_B-\Pi_A}_{\Ls(H)}
    \leq
    C_A|A\triangle B|.
\end{equation}
\end{theorem}

\begin{proof}
We set for convenience $D:=A\triangle B$ and $\eta:=\chi_B-\chi_A$. We write the Riccati equation for \(B\) as a perturbation of the Riccati
equation for \(A\). If $\Delta\Pi:=\Pi_B-\Pi_A$, then due to \cref{lem:riccati_perturbation_identity} \(\Delta\Pi\) satisfies
\[
\begin{aligned}
    &
    (\A_A^{\cl})^*\Delta\Pi
    +
    \Delta\Pi\A_A^{\cl}
    \\
    &=
    \frac1\alpha
    \Big[
        \Pi_A\eta\Pi_A
        +
        \Pi_A\eta\Delta\Pi
        +
        \Delta\Pi\eta\Pi_A
        +
        \Delta\Pi\chi_A\Delta\Pi
        +
        \Delta\Pi\eta\Delta\Pi
    \Big].
\end{aligned}
\]
We prove that this equation has a unique small solution by a Banach fixed-point
argument. Throughout this section, \(C>0\) denotes a generic constant that may change
from line to line. The constant may depend on the fixed reference support
\(A\) and on the structural constants of the problem, but it is
independent of the perturbing support \(B\), of \(D=A\triangle B\), and of the
auxiliary variables being estimated.

Since Riccati operators are self-adjoint, we work on the closed
subspace of self-adjoint bounded operators on \(H\). Define the map \(\Phi_B\) by
\begin{equation} \label{eq_def_Phi_B}
    \Phi_B(X)
    :=
    \LA^{-1}
    \left[
        \frac1\alpha
        \Big(
            \Pi_A\eta\Pi_A
            +
            \Pi_A\eta X
            +
            X\eta\Pi_A
            +
            X\chi_A X
            +
            X\eta X
        \Big)
    \right],
\end{equation}
where  $\LA X=(\A_A^{\cl})^*X+X\A_A^{\cl}$.

\paragraph{Well-posedness of \(\Phi_B\)}
For self-adjoint $X$ the expression inside the brackets in \eqref{eq_def_Phi_B} is self-adjoint.
%since \(\Pi_A\eta X\) and \(X\eta \Pi_A\) are adjoint to each other and the remaining terms are self-adjoint.
Hence, by the Lyapunov representation in \cref{lem:standard_lyapunov_inverse} also \(\Phi_B(X)\) is self-adjoint.  Next, we  estimate the terms in \(\Ls(H)\).
First, we have
\begin{equation}
\label{eq:thm15_1}
\norm{\Pi_A\eta\Pi_A}_{\Ls(H)}
 =\sup_{\norm{v}_H\leq1}\left|\ip{\Pi_Av}{\eta\Pi_Av}_H\right|
 \leq |D|\norm{\Pi_A}_{\Ls(H,L^\infty(\Omega))}^2\leq C|D|.
\end{equation}
Since $\Pi_A\eta X + X\eta\Pi_A$ is self-adjoint, we can estimate
\begin{equation}
\label{eq:thm15_2}
\begin{aligned}
 \norm{\Pi_A\eta X+ X\eta\Pi_A}_{\Ls(H)} & = 2\sup_{\norm{v}_H \le 1} \left|\ip{\Pi_A v}{\eta Xv}_H\right| \\
 &\le 2\sup_{\norm{v}_H \le 1} \norm{\Pi_Av}_{L^\infty(\Omega)} |D|^{1/2} \norm{Xv}_H
 \le C |D|^{1/2} \norm{X}_{\Ls(H)}.
 \end{aligned}
\end{equation}
As multiplication by characteristic functions is bounded on
\(H\) with norm at most one, we have
\begin{equation}
\label{eq:thm15_4}
\norm{X\chi_A X + X\eta X}_{\Ls(H)} = \norm{X\chi_B X }_{\Ls(H)}  \leq \norm{X}_{\Ls(H)}^2.
\end{equation}
By \cref{lem:standard_lyapunov_inverse} and using \eqref{eq:thm15_1}-\eqref{eq:thm15_4}, we get the bound
\begin{equation}
\label{eq:bound_of_Phi_B}
\begin{aligned}
    \norm{\Phi_B(X)}_{\Ls(H)}
    &\leq
    C\Big[
        |D|^{1/2}
        +
        \norm{X}_{\Ls(H)}
    \Big]^2.
    \end{aligned}
 \end{equation}
Take $R>C$.
We now show that \(\Phi_B\) is a contraction on the ball
\[
    \mathcal B_D
    :=
    \left\{
        X=X^*\in\Ls(H):
        \norm{X}_{\Ls(H)}
        \leq
        R|D|
    \right\}.
\]
Take $X \in \mathcal B_D$. Then by the above estimates
\[
    \norm{\Phi_B(X)}_{\Ls(H)}
    \leq
    C |D| ( 1 + |D|^{1/2} R)^2.
\]
Hence, for $|D|$ sufficiently small, we get $\norm{\Phi_B(X)}_{\Ls(H)} \le R |D|$.
Thus \(\Phi_B(\mathcal B_D)\subset\mathcal B_D\).

\paragraph{Contraction} Let \(X,Y\in\mathcal B_D\) be arbitrary. Estimating the difference
\(\Phi_B(X)-\Phi_B(Y)\) in the same way as in \eqref{eq:thm15_2}-\eqref{eq:thm15_4}  gives us
\begin{equation}
\label{eq:Lipschitz_Phi_B}
\begin{aligned}
    \norm{\Phi_B(X)-\Phi_B(Y)}_{\Ls(H)}
    &\leq
    C\Big[
        |D|^{1/2}
        +
        \norm{X}_{\Ls(H)}
        +
        \norm{Y}_{\Ls(H)}
    \Big]
    \norm{X-Y}_{\Ls(H)}.
\end{aligned}
\end{equation}
Since \(X,Y\in\mathcal B_D\), we also have $\norm{X}_{\Ls(H)}+\norm{Y}_{\Ls(H)}\leq 2R|D|$. Hence the Lipschitz constant $C\Big[|D|^{1/2} +\norm{X}_{\Ls(H)}+\norm{Y}_{\Ls(H)} \Big]$ tends to zero as \(|D|\to0\). Thus, after
possibly decreasing \(|D|\), \(\Phi_B\) is a contraction on
\(\mathcal B_D\). The Banach fixed-point theorem gives a unique fixed point $\Delta\Pi_B\in\mathcal B_D$ of $\Phi_B$.
In particular, $\norm{\Delta\Pi_B}_{\Ls(H)}\leq R|D|=R|A\triangle B|$.

\paragraph{Identification with the stabilizing Riccati solution}
It remains to show that the fixed point constructed above is the stabilizing
Riccati solution corresponding to the support \(B\). Set \(\widetilde\Pi_B:=\Pi_A+\Delta\Pi_B\). By construction,
\(\widetilde\Pi_B\) solves the Riccati equation for the support \(B\). Its
closed-loop generator is
\[
    \A-\frac1\alpha\chi_B\widetilde\Pi_B
    =
    \A_A^{\cl}
    -
    \frac1\alpha
    \left[
        (\chi_B-\chi_A)\Pi_A+\chi_B\Delta\Pi_B
    \right].
\]
The perturbation is small in \(\Ls(H)\) since by H\"older inequality and
the construction of $\Delta\Pi_B$ we get
%due to \cref{lem:localization_symmetric_difference} and \eqref{eq:Pi_difference_rate}, we have
\[
\begin{aligned}
    \norm{(\chi_B-\chi_A)\Pi_A}_{\Ls(H)}
    &\leq  |A\triangle B|^{1/2} \norm{\Pi_A}_{\Ls(H,L^\infty(\Omega))},\\
    \norm{\chi_B\Delta\Pi_B}_{\Ls(H)}
    &\leq R|A\triangle B|.
\end{aligned}
\]
Hence, by the same bounded-perturbation argument as in the proof of
\cref{lem:local_stabilizability}, this closed loop is exponentially stable
for sufficiently small \(|A\triangle B|\).

Now, we define the feedback associated with \(\widetilde\Pi_B\) by  $\widetilde\KB := -\frac1\alpha\chi_B\widetilde\Pi_B$, and let
\[
    \widetilde\A_B^{\cl}
    :=
    \A+\chi_B\widetilde\KB
    =
    \A-\frac1\alpha\chi_B\widetilde\Pi_B
\]
be the corresponding closed-loop generator. We denote by
\(\widetilde S_B^{\cl}(t)\) the semigroup generated by
\(\widetilde\A_B^{\cl}\). Since this semigroup is exponentially stable, the
closed-loop Lyapunov representation gives
\[
    \widetilde\Pi_B
    =
    \int_0^\infty
    \widetilde S_B^{\cl}(t)^*
    \left(
        Q^*Q+\alpha\widetilde\KB^*\widetilde\KB
    \right)
    \widetilde S_B^{\cl}(t)
    \,\dd t .
\]
The integrand is nonnegative, and hence \(\widetilde\Pi_B\geq0\). Therefore
\(\widetilde\Pi_B\) is a nonnegative stabilizing solution of the Riccati
equation for the support \(B\). By uniqueness of the nonnegative stabilizing solution of the algebraic
Riccati equation, see, e.g.,
\cite[Thm.~6.2.7]{CurtainZwart1995} and \cite[Pt.~V, Ch.~1, Secs.~2--3]{BensoussanDaPratoDelfourMitter2007}, we obtain \(\widetilde\Pi_B=\Pi_B\). Consequently, $\Pi_B- \Pi_A=\Delta\Pi_B$, and therefore, $\norm{\Pi_B-\Pi_A}_{\Ls(H)}\leq R|A\triangle B|$ holds.
% Finally, using \eqref{eq:psi_difference_rate}, we can write
% \[
% \begin{aligned}
%     \norm{(\Pi_B-
%     \Pi_A)y_A}_{L^2(0,\infty;H)}
%     &\leq
%     \norm{\Pi_B-
%     \Pi_A}_{\Ls(H)}
%     \norm{y_A}_{L^2(0,\infty;H)}
%     \leq
%     C|A\triangle B| \norm{y_A}_{L^2(0,\infty;H)},
% \end{aligned}
% \]
% and the proof is finished.
\end{proof}

\begin{corollary}[Remainder estimate]
\label{cor:remainder_estimate}
Let \(A\in\Aad\)
with \(\Pi_A\in\Ls(H,L^\infty(\Omega))\).
Then there
exist constants \(C_{\rm rem,A}>0\) and \(\delta_A>0\) such that for every stabilizing
support \(B\in\Aad\) with $|A\triangle B|<\delta_A$, the following
estimate of the first-order expansion of $J$ is satisfied:
\begin{equation}
\label{eq:first_order_expansion}
   \left| J(B)-J(A) -
    \int_\Omega
    (\chi_B-\chi_A)(x)G_A(x)\,\dd x
    \right|
    \le
    C_{\rm rem,A} |A\triangle B|^{3/2}.
\end{equation}
\end{corollary}

\begin{proof}
Set \(D:=A\triangle B\). Due to the assumptions,  we have
\(\Pi_A\in\Ls(H,L^\infty(\Omega))\). Hence
\(\psi_A:=\Pi_Ay_A\in L^2(0,\infty;L^\infty(\Omega))\). Moreover,
\(\psi_B-\psi_A=(\Pi_B-\Pi_A)y_A\), and \cref{thm:local_branch_and_sensitivity}
implies
\[
    \norm{\psi_B-\psi_A}_{L^2(0,\infty;H)}
    \le C |D|,
\]
with $C$ depending only on $A$.
Since \(\eta_{A,B}:=\chi_B-\chi_A\) is supported on \(D\), we also have
\[
    \norm{\eta_{A,B}\psi_A}_{L^2(0,\infty;H)}
    \leq
    |D|^{1/2}
    \norm{\psi_A}_{L^2(0,\infty;L^\infty(\Omega))}
    =
    C|D|^{1/2}.
\]
Using \eqref{eq:remainder_explicit},
% Cauchy's inequality, and
% \(\norm{\chi_B(\psi_B-\psi_A)}_{L^2(0,\infty;H)}
% \leq \norm{\psi_B-\psi_A}_{L^2(0,\infty;H)}\),
we obtain
\[
\begin{aligned}
    |\mathcal R_A(B)|
    &\leq
    C
    \norm{\eta_{A,B}\psi_A}_{L^2(0,\infty;H)}
    \norm{\psi_B-\psi_A}_{L^2(0,\infty;H)}
    +
    C
    \norm{\psi_B-\psi_A}_{L^2(0,\infty;H)}^2
    \\
    &\leq
    C|D|^{3/2}
    +
    C|D|^{2}.
\end{aligned}
\]
Thus \(\mathcal R_A(B)=O(|D|^{3/2})=O(|A\triangle B|^{3/2})\). Combining this with
\Cref{lem:exact_expansion_remainder} gives \eqref{eq:first_order_expansion}.
\end{proof}

\section{Topological derivative and descent algorithm}
\label{sec:topological_algorithm}

We now use the exact expansion and the remainder estimate to identify the
topological derivative of the support functional
for the class \(\Aad\). The resulting switching
function is then written in a computable Lyapunov form and used to define a
symmetric-difference descent algorithm.

\subsection{Topological derivative}
\label{subsec:topological_derivative}

For \(x\in\Omega\) and \(\varepsilon>0\), we write  $B_\varepsilon(x):=\{\xi\in\Omega:\ |\xi-x|<\varepsilon\}$. Thus \(B_\varepsilon(x)\) denotes the ball of radius \(\varepsilon\) centered
at \(x\), understood relative to the spatial domain \(\Omega\).

\begin{theorem}[Topological derivative]
\label{thm:topological_derivative}
Let \(A\in\Aad\) be a stabilizing support
such that \(\Pi_A\in\Ls(H,L^\infty(\Omega))\).
Then the following assertions hold
for almost every \(x\in\Omega\):
\begin{enumerate}[label=\textup{(\roman*)}]
    \item If \(x\in\Omega\setminus A\), then
    \[
        \lim_{\varepsilon\to0}
        \frac{J(A\cup B_\varepsilon(x))-J(A)}
        {|B_\varepsilon(x)|}
        =
        G_A(x).
    \]
%     Hence \(D_T^{\mathrm{add}}J(A)(x)=G_A(x)\). <- was eh undefined

    \item If \(x\in A\), then
    \[
        \lim_{\varepsilon\to0}
        \frac{J(A\setminus B_\varepsilon(x))-J(A)}
        {|B_\varepsilon(x)|}
        =
        -G_A(x).
    \]
%     Hence \(D_T^{\mathrm{rem}}J(A)(x)=-G_A(x)\).
\end{enumerate}
\end{theorem}
\begin{proof}
By the Lebesgue differentiation theorem,
there exists a null set \(N\subset\Omega\) such that every
\(x\in\Omega\setminus N\) is a Lebesgue point of \(G_A\) and $\chi_AG_A$.
We fix such an \(x\). By \cref{lem:local_stabilizability}, the
supports \(A\cup B_\varepsilon(x)\) and \(A\setminus B_\varepsilon(x)\)
are stabilizable for sufficiently small \(\varepsilon\).
First let \(x\in\Omega\setminus A\). Then we have
\[
    (A\cup B_\varepsilon(x))\triangle A
    =
    B_\varepsilon(x)\setminus A, \quad \text{ and }
    \quad
    \chi_{A\cup B_\varepsilon(x)}-\chi_A
    =
    \chi_{B_\varepsilon(x)\setminus A}
    =
    \chi_{B_\varepsilon(x)}(1-\chi_A).
\]
Hence, by \eqref{eq:first_order_expansion}, we obtain
\[
    J(A\cup B_\varepsilon(x))-J(A)
    =
    \int_{B_\varepsilon(x)}  (1-\chi_A) G_A(\xi)\,\dd\xi
    +
    o(|B_\varepsilon(x)|).
\]
Dividing by \(|B_\varepsilon(x)|\) and using the Lebesgue differentiation
theorem yields
\[
    \lim_{\varepsilon\to0}
    \frac{J(A\cup B_\varepsilon(x))-J(A)}
    {|B_\varepsilon(x)|}
    =
    (1-\chi_A(x)) G_A(x) = G_A(x).
\]
The removal case is analogous.
If \(x\in A\), then we have
\[
    (A\setminus B_\varepsilon(x))\triangle A
    =
    A\cap B_\varepsilon(x),
    \qquad
    \chi_{A\setminus B_\varepsilon(x)}-\chi_A
    =
    -\chi_{A\cap B_\varepsilon(x)}
    = -\chi_A \chi_{B_\varepsilon(x)}.
\]
The first-order expansion gives
\[
    J(A\setminus B_\varepsilon(x))-J(A)
    =
    -\int_{B_\varepsilon(x)}\chi_A G_A(\xi)\,\dd\xi
    +
    o(|B_\varepsilon(x)|).
\]
Dividing by \(|B_\varepsilon(x)|\) again gives
\[
    \lim_{\varepsilon\to0}
    \frac{J(A\setminus B_\varepsilon(x))-J(A)}
    {|B_\varepsilon(x)|}
    =
    -G_A(x).
\]
This proves the result.
\end{proof}

The theorem shows that actuator material should be added where \(G_A<0\) and
removed where \(G_A>0\). Equivalently, the sign of \(G_A\) gives the
first-order effect of changing the support locally.

\begin{theorem}[Pontryagin maximum principle]
\label{thm:pmp}
Let \(A\in\Aad\) be a stabilizing support
with \(\Pi_A\in\Ls(H,L^\infty(\Omega))\),
so that $\chi_A$ is locally optimal in $L^1(\Omega)$, i.e.,
there is $\rho>0$ such that
\[
 J(A) \le J(B) \quad \forall B\in \Aad: \  |A\triangle B|=\norm{\chi_A - \chi_B}_{L^1(\Omega)} <\rho.
\]
Then it holds
\begin{equation}\label{eq_optimality_condition}
 (1-2 \chi_A) G_A \ge 0 \quad \text{ a.e.\@ in }\Omega.
\end{equation}
\end{theorem}
\begin{proof}
Take $x\in A$ such that the result of \cref{thm:topological_derivative} is valid.
Since $\chi_A$ is locally optimal,
it holds $J(A\setminus B_\varepsilon(x))-J(A) \ge 0$ for all $\varepsilon>0$ sufficiently small,
which yields $G_A(x) \le0$.
Similarly, $x\in \Omega\setminus A$ implies $G_A(x) \ge 0$.
Consequently, $(1-2 \chi_A(x)) G_A \ge 0$.
\end{proof}

The optimality condition \eqref{eq_optimality_condition} of \cref{thm:pmp} can be written in the maximum form
\[
 \chi_A(x) \in \operatorname{argmax}_{v\in \{0,1\} } \left( {-}v\cdot G_A(x) \right) \quad \text{ for a.a.\@ } x \in \Omega.
\]
The argmax contains both \(0\) and \(1\) when \(G_A(x)=0\).

\subsection{Lyapunov-based switching function}

In computations, %For numerical purposes,
the infinite-horizon term in the switching function is
not evaluated by direct time integration. Instead, it is computed through a
closed-loop Lyapunov equation. Let
\(\A_A^{\cl}=\A-\alpha^{-1}\chi_A\Pi_A\) be the closed-loop generator and let
\(S_A^{\cl}(t)=e^{t\A_A^{\cl}}\). For a given initial datum \(y_0\), define
the trajectory Gramian \(Z_A\) as the solution of
\begin{equation}
\label{eq:ZA_lyapunov_algorithm}
    \A_A^{\cl}Z_A+Z_A(\A_A^{\cl})^*+y_0\otimes y_0=0.
\end{equation}
The equation is understood weakly, by testing against elements of
\(\Dom((\A_A^{\cl})^*)\).  Here \(y_0\otimes y_0\) denotes the rank-one
operator on \(H\) given by
\((y_0\otimes y_0)z=\ip{z}{y_0}_H y_0\) for all \(z\in H\). Since the closed-loop
semigroup is exponentially stable, \eqref{eq:ZA_lyapunov_algorithm} is
equivalently represented by
\[
    Z_A
    =
    \int_0^\infty
    S_A^{\cl}(t)y_0\otimes S_A^{\cl}(t)y_0\,\dd t
    =
    \int_0^\infty
    y_A(t)\otimes y_A(t)\,\dd t .
\]
Thus \(Z_A\) collects the time-correlation of the closed-loop trajectory
\(y_A(t)=S_A^{\cl}(t)y_0\). Applying \(\Pi_A\) from both sides gives
\[
    \Pi_AZ_A\Pi_A
    =
    \int_0^\infty
    (\Pi_Ay_A(t))\otimes(\Pi_Ay_A(t))\,\dd t,
\]
where we used the fact that \(\Pi_A\) is self-adjoint.
This positive operator is trace class, since
\(\int_0^\infty\norm{\Pi_Ay_A(t)}_H^2\,\dd t<\infty\).
It follows that
\begin{equation} \label{eq_trace_identity}
 \int_\Omega \int_0^\infty
    |\Pi_Ay_A(t)(x)|^2\,\dd t \, \dd x = \operatorname{trace}( \Pi_AZ_A\Pi_A ).
\end{equation}
Let us assume $\Pi_A \in \Ls(H,L^\infty(\Omega))$. Then for almost all $x\in \Omega$, we can write
\[
 \int_0^\infty |\Pi_Ay_A(t)(x)|^2\,\dd t
 =\lim_{\epsilon\searrow0}\int_0^\infty
 \left(\frac1{|B_\epsilon(x)|}\int_{B_\epsilon(x)}
 \Pi_Ay_A(t)(\xi)\,\dd\xi\right)^2\dd t,
\]
which follows from the Lebesgue differentiation theorem, Fubini's theorem, and the dominated convergence theorem. Indeed, using the exponential decay from \eqref{eq:closed_loop_exp_stability} the averages can be  bounded by
$  \norm{\Pi_A}_{\Ls(H,L^\infty(\Omega))}  M_Ae^{-\gamma_A t}\norm{y_0}_H$.
Consequently the limit
\[
 \lim_{\epsilon\searrow0}  \frac1{|B_\epsilon(x)|^2}  \ip{ \chi_{B_\epsilon(x)} }{ \Pi_AZ_A\Pi_A\chi_{B_\epsilon(x)} }_H
\]
exists for almost all $x\in \Omega$. Motivated by \eqref{eq_trace_identity}, we call the resulting function the spatial diagonal of $\Pi_AZ_A\Pi_A$, i.e.,
\[
 \diag(\Pi_AZ_A\Pi_A)(x) := \lim_{\epsilon\searrow0}  \frac1{|B_\epsilon(x)|^2}  \ip{ \chi_{B_\epsilon(x)} }{ \Pi_AZ_A\Pi_A\chi_{B_\epsilon(x)} }_H
%  =  \int_0^\infty  |\Pi_Ay_A(t)(x)|^2\,\dd t
\]
for almost all $x\in \Omega$.
Hence, the spatial diagonal of
\(\Pi_AZ_A\Pi_A\) represents the pointwise infinite-horizon quantity
\[
    \diag(\Pi_AZ_A\Pi_A)(x)
    =
    \int_0^\infty
    |\Pi_Ay_A(t)(x)|^2\,\dd t .
\]
This identity is the infinite-dimensional analog of the matrix fact that the
diagonal of \(vv^\top\) consists of the componentwise squares of \(v\). Consequently, the switching function can be computed as
\begin{equation}
\label{eq:GA_computable}
    G_A(x)
    =
    \beta(x)
    -
    \frac1\alpha
    \diag(\Pi_AZ_A\Pi_A)(x).
\end{equation}
After spatial discretization, \(Z_A\) becomes the solution of a matrix
Lyapunov equation and \(\diag(\Pi_AZ_A\Pi_A)\) becomes the vector of diagonal
entries of the corresponding matrix product, up to the mass-matrix convention
of the chosen finite-element discretization.

\subsection{Descent update}
\label{subsec:descent_algorithm}

To obtain simultaneous additions and removals, let \(D\subset\Omega\) be a
measurable perturbation set. Since \(\chi_{A\triangle D}=|\chi_A-\chi_D|
=\chi_A+\chi_D-2\chi_A\chi_D\), we have $\chi_{A\triangle D}-\chi_A=(1-2\chi_A)\chi_D$. Thus, the first-order expansion of \cref{cor:remainder_estimate}
gives
\begin{equation}
\label{eq:expansion_with_triangle}
    J(A\triangle D)-J(A)
    =
    \int_D(1-2\chi_A(x))G_A(x)\,\dd x
    +
    o(|D|).
\end{equation}
We therefore define the flip indicator
\begin{equation}
\label{eq:flip_indicator_algorithm}
    \Theta_A(x):=(1-2\chi_A(x))G_A(x).
    \end{equation}
The condition $\Theta_A(x)\geq0$ for a.e.\@ $x \in\Omega$ is the natural first-order optimality condition for the support variable of \cref{thm:pmp}.
% It is equivalent to \(G_A\geq0\) a.e.\@ on \(\Omega\setminus A\) and \(G_A\leq0\) a.e.\@ on \(A\).
Thus, if \(\Theta_A<0\) on a set of positive measure $D$, flipping part of this region, i.e., replacing $A$ by $A\triangle D$, gives a first-order descent direction.
We define
\begin{equation}
\label{eq:stationarity_defect}
    q_A:=\min\{\Theta_A,0\},\qquad
    \eta_A:=\norm{q_A}_{L^1(\Omega)},\qquad
    R_A:=\{x\in\Omega:q_A(x)<0\}.
\end{equation}
Then the first-order condition can be expressed equivalently as \(\eta_A=0\).

We use an Armijo line-search in the measure of the changed set, following
\cite[Sec.~6]{Wachsmuth2024}. At a nonstationary iterate \(A_k\), abbreviate
\(q_k=q_{A_k}\), $\Theta_k = \Theta_{A_k}$, \(\eta_k=\eta_{A_k}>0\), and \(R_k=R_{A_k}\).
For a trial fraction \(s\in(0,1]\), choose a measurable set satisfying
\begin{equation}
\label{eq:measure_controlled_selection}
    D_{k,s}\subset R_k,\qquad
    |D_{k,s}|\le s|R_k|,\qquad
    \int_{D_{k,s}}q_k\dd x\leq s\int_{R_k}q_k\dd x =  -s\eta_k .
\end{equation}
Since the Lebesgue measure is nonatomic, there exists a set $D_{k,s}$ satisfying \eqref{eq:measure_controlled_selection}
with equality owing to the celebrated Lyapunov convexity theorem.

We denote the trial support by \(A_{k,s}:=A_k\triangle D_{k,s}\). The value of $s$ as well as $A_{k,s}$ are accepted, i.e., $A_{k+1} := A_{k,s}$,  if the
Armijo-like descent condition
\begin{equation}
\label{eq:accepted_armijo_condition}
    A_{k,s}\in\Aad,\qquad
    J(A_{k,s})\leq J(A_k)+\sigma\int_{D_{k,s}}q_k\dd x
\end{equation}
is satisfied, where $\sigma\in (0,1)$ is a fixed parameter.
As in \cite[Lem.~6.5]{Wachsmuth2024}, we can show that there exist values of $s$ such that
this condition is satisfied.

\begin{lemma}
\label{lem:backtracking_works}
Let $A_k \in \Aad$ be given with \(\Pi_{A_k}\in\Ls(H,L^\infty(\Omega))\) and with $\eta_k>0$, i.e., it does not satisfy the necessary optimality condition. Then there is $s_0\in (0,1)$ such that for all $s\in (0,s_0)$ and all $D_{k,s}$ satisfying \eqref{eq:measure_controlled_selection}
the set $A_{k,s}=A_k\triangle D_{k,s}$ satisfies the descent condition \eqref{eq:accepted_armijo_condition}.
\end{lemma}
\begin{proof}
By construction $A_{k,s} \triangle A_k = D_{k,s}$.
Due to \cref{lem:local_stabilizability}, $A_{k,s}$ belongs to $\Aad$ provided that $|D_{k,s}|$ is small enough.
In addition, the first-order expansion of \cref{cor:remainder_estimate} is valid if $|D_{k,s}|$ is small enough.
Consequently,
there is $s_0>0$ such that $A_{k,s}$ belongs to $\Aad$ and $D_{k,s}$ satisfies \eqref{eq:expansion_with_triangle} for all $s\in (0,s_0)$.

By construction, $q_k =\Theta_k<0$ on $R_k\supset D_{k,s}$.
Hence we get with \eqref{eq:expansion_with_triangle}, \eqref{eq:flip_indicator_algorithm}, \eqref{eq:measure_controlled_selection}
% \begin{multline*}
\[
 \begin{aligned}
 J(A_k\triangle D_{k,s})-J(A_k) - \sigma\int_{D_{k,s}}q_k\dd x
& \le (1-\sigma) \int_{D_{k,s}} \Theta_k(x)\dd x
    +
    o(|D_{k,s}|)\\
& \le -(1-\sigma) s \eta_k + o(s).
    \end{aligned}\\
\]%\end{multline*}
Since $\eta_k>0$, the right-hand side is negative for small enough $s$.
\end{proof}

\Cref{lem:backtracking_works} ensures that the inner backtracking loop of
\cref{alg:lyapunov_support_descent} terminates at each outer iteration.
Each iteration of the inner loop requires one Riccati solve to compute the value function at the trial support.
In addition, one Lyapunov solve is performed in each outer iteration to compute the derivative $G_A$.
We remark that no truncation of an
infinite-horizon time integral is needed.

\begin{algorithm}[htbp]
\caption{Lyapunov-based support descent with measure-controlled Armijo search}
\label{alg:lyapunov_support_descent}
\begin{algorithmic}[1]
\State Choose \(A_0\in\Aad\), \(0<\rho,\sigma<1\), and
\(\varepsilon_{\rm stop}\geq0\).
\For{\(k=0,1,2,\ldots\)}
    \State Solve the Riccati equation for \(\Pi_{A_k}\) and form
    \(\A_{A_k}^{\cl}=\A-\alpha^{-1}\chi_{A_k}\Pi_{A_k}\).
    \State Solve \eqref{eq:ZA_lyapunov_algorithm} for \(Z_k\).
    \State Compute \(G_k,\Theta_k,q_k,\eta_k,R_k\) using
    \cref{eq:GA_computable,eq:flip_indicator_algorithm,eq:stationarity_defect}.
    \State \textbf{if} \(\eta_k\leq\varepsilon_{\rm stop}\)
    \textbf{then} \Return \(A_k,\eta_k\).
    \For{\(s=1,\rho,\rho^2,\ldots\)}
        \State Choose \(D_{k,s}\) satisfying
        \eqref{eq:measure_controlled_selection};
        set \(A_{k,s}=A_k\triangle D_{k,s}\).
        \If{the acceptance condition
            \eqref{eq:accepted_armijo_condition} holds}
            \State \(s_k\gets s\),
            \(A_{k+1}\gets A_{k,s_k}\); \textbf{break}.
        \EndIf
    \EndFor
\EndFor
\end{algorithmic}
\end{algorithm}
The termination criterion is motivated by the following result, which is a consequence of the one-sided estimate of \cref{prop:value_one_sided_estimate},
see also \cite[Cor.~4.5]{Wachsmuth2024}.

\begin{lemma}
\label{lem:computable_error_estimate_V}
For \(A\in\Aad\), we have \(0\leq J(A)-\inf_{B\in\Aad}J(B)\leq\eta_A\).
\end{lemma}
\begin{proof}
Let \(B\in\Aad\). By \eqref{eq:switching_function}, the right-hand side of
\eqref{eq:value_one_sided_estimate} is
\(\int_\Omega(\chi_B-\chi_A)(G_A-\beta)\,\dd x\). Hence
\[
 J(A)-J(B)
 =V_A(y_0)-V_B(y_0)-\int_\Omega(\chi_B-\chi_A)\beta\,\dd x
 \leq-\int_\Omega(\chi_B-\chi_A)G_A\,\dd x.
\]
Since \(\chi_B-\chi_A=(1-2\chi_A)\chi_{A\triangle B}\),
\eqref{eq:flip_indicator_algorithm} and \eqref{eq:stationarity_defect} give
\[
 (\chi_B-\chi_A)G_A=\chi_{A\triangle B}\Theta_A
 \geq\min\{\Theta_A,0\}=q_A.
\]
Therefore \(J(A)-J(B)\leq-\int_\Omega q_A\,\dd x=\eta_A\).
Taking the supremum of \(J(A)-J(B)\) over \(B\in\Aad\) gives the upper
bound. The lower bound follows because \(A\in\Aad\).
\end{proof}

In the following, we assume that \cref{alg:lyapunov_support_descent} generates an infinite sequence of iterates.

\begin{lemma}
\label{lem:basic_convergence_estimate}
 Let $(A_k)$ be an infinite sequence of iterates of \cref{alg:lyapunov_support_descent}.
 Then it holds
 \[
  \sum_{k=1}^\infty s_k \eta_k < +\infty.
 \]
\end{lemma}
\begin{proof}
 The proof is similar to that of \cite[Lem.~6.6]{Wachsmuth2024}.
 From \eqref{eq:measure_controlled_selection} and \eqref{eq:accepted_armijo_condition} we get
\[
 J(A_{k+1}) - J(A_k)  \le  \sigma \int_{ D_{k,s_k}}q_k \dd x \le - \sigma s_k \eta_k.
\]
Since $J$ is bounded from below, we can sum this inequality for $k=0, \dots, \infty$ to prove the claim.
\end{proof}

To prove the convergence of the residuals $\eta_k$ to zero, one needs additional assumptions. Such a result was proven in \cite[Thm.~6.7]{Wachsmuth2024}, where a uniform remainder estimate was used. In our situation, this amounts to the validity of a uniform bound on the remainder in the first-order expansion of $J$ in \cref{cor:remainder_estimate}.

\begin{theorem}
\label{thm_convergence_under_uniform_remainder}
Let $(A_k)$ be an infinite sequence of iterates of \cref{alg:lyapunov_support_descent}.
Suppose that there are
constants \(R>0\) and \(\delta>0\) such that,
for every \(k\) and every measurable \(B\subset\Omega\) with
\(|A_k\triangle B|<\delta\),
the following two conditions are satisfied:
\begin{enumerate}
 \item \(B\) is exponentially stabilizing and hence belongs to \(\Aad\),
 \item
$%\begin{equation}
%\label{eq:first_order_expansion_uniform}
   \left| J(B)-J(A_k) -
    \int_\Omega
    (\chi_B-\chi_{A_k})(x)G_{A_k}(x)\,\dd x
    \right|
    \le
    R |A_k\triangle B|^{3/2}.
$%\end{equation}
\end{enumerate}
Then it holds $\eta_k \to 0$, i.e., the residual in the optimality condition converges to zero, and $(A_k)$ is a minimizing sequence.
\end{theorem}
\begin{proof}
Set $\tilde s:=\min\{1,\rho\delta/|\Omega|\}$.
 Take $k$ such that $s_k< \tilde s$. Due to the backtracking strategy, the step $s := s_k/\rho$
 violates the condition \eqref{eq:accepted_armijo_condition}, i.e.,
\[
0< J(A_k\triangle D_{k,s})-J(A_k) - \sigma\int_{D_{k,s}}q_k\dd x.
\]
By \eqref{eq:measure_controlled_selection} and the choice of $\tilde s$, it follows $|D_{k,s}| \le s |\Omega|< \delta$.
Hence, we can apply the uniform error estimate from the assumptions
and derive similarly to  \eqref{eq:expansion_with_triangle},
\[\begin{aligned}
  J(A_k\triangle D_{k,s})-J(A_k)
&\le
    \int_{D_{k,s}} (1-2\chi_{A_k}(x))G_{A_k} (x)\,\dd x
    +
    R|D_{k,s}|^{3/2}\\
    &=\int_{D_{k,s}} \Theta_{A_k}(x)\,\dd x
    +
    R|D_{k,s}|^{3/2}.
\end{aligned}\]
Since $q_k = \Theta_k <0 $ on $D_{k,s}$, cf. \eqref{eq:stationarity_defect},
we obtain the inequality
\[
 0<  (1-\sigma)\int_{D_{k,s}}q_k\dd x +  R|D_{k,s}|^{3/2}
\]
By the conditions \eqref{eq:measure_controlled_selection} on $D_{k,s}$, we get
\[
 0< (1-\sigma)\int_{D_{k,s}}q_k\dd x +  R|D_{k,s}|^{3/2}
 \le -(1-\sigma)s\eta_k+Rs^{3/2}|\Omega|^{3/2},
\]
which yields the bound
\[
s_k=\rho s\geq\rho\left(\frac{1-\sigma}{R|\Omega|^{3/2}}\right)^2\eta_k^2.
\]
Using this lower bound in \cref{lem:basic_convergence_estimate} results in
$
\sum_{k:\, s_k\ge \tilde s} \eta_k + \sum_{k:\, s_k<\tilde s} \eta_k^3 < +\infty
$, so that $\eta_k\to 0$ follows.

By \cref{lem:computable_error_estimate_V},
\(J(A_k)\to\inf_{B\in\Aad}J(B)\). Thus \((A_k)\) is a minimizing sequence
for the total objective.
\end{proof}

Let us provide a sufficient condition for the uniform remainder estimate that was used in the previous theorem.

\begin{theorem}
\label{thm_sufficient_condition_for_under_uniform_remainder}
Let $(A_k)$ be an infinite sequence of iterates of \cref{alg:lyapunov_support_descent}.
Assume that the following conditions are satisfied:
\begin{enumerate}
 \item $(M_{A_k})$ is bounded from above, and $(\gamma_{A_k})$ is bounded from below by a positive number, which are the sequences of the constants in the exponential closed-loop estimate \eqref{eq:closed_loop_exp_stability},
\item $( \Pi_{A_k} )$ is uniformly bounded in $\Ls(H,L^\infty(\Omega))$.
\end{enumerate}
Then the assumptions of \cref{thm_convergence_under_uniform_remainder} are satisfied.
\end{theorem}
\begin{proof}
Due to \cref{lem:closed_loop_W_regular}, $(y_{A_k})$ is bounded in $W(0,+\infty)$.
An inspection of the proof of \cref{lem:standard_lyapunov_inverse} shows that $(\mathfrak L_{A_k}^{-1})$ is uniformly bounded in $\Ls(\Ls(H),\Ls(H))$.
By \cref{lem:local_stabilizability}, we get a uniform constant $\delta>0$ such that every set $B$ with $|A_k \triangle B| <\delta$ for some $k$ is a stabilizing support.

Now let us briefly inspect the proof of \cref{thm:local_branch_and_sensitivity}:
the constants $C$ in
\eqref{eq:bound_of_Phi_B} and in  \eqref{eq:Lipschitz_Phi_B} depend on $\norm{\mathfrak L_{A_k}^{-1}}_{\Ls(\Ls(H),\Ls(H))}$ and $\norm{\Pi_{A_k}}_{\Ls(H,L^\infty(\Omega))}$,
and hence are uniformly bounded for all $A_k$.
Following the contraction argument, one finds that the sequences of constants $(C_{A_k})$ and $(\delta_{A_k})$ provided by \cref{thm:local_branch_and_sensitivity}
are bounded uniformly from above and below by a positive number, respectively.
The proof of \cref{cor:remainder_estimate} now shows that the sequences $(C_{\rm rem,A_k})$ and $(\delta_{A_k})$ provided by this result are bounded
from above and below by a positive number, respectively.
Hence, the assumptions of \cref{thm_convergence_under_uniform_remainder} are satisfied.
\end{proof}

\subsection{Uniform bounds}
\label{subsec:uniform_bounds_stationarity}

We first establish bounds that are uniform with respect to the actuator, so that \cref{thm_sufficient_condition_for_under_uniform_remainder,thm_convergence_under_uniform_remainder}
apply. For this purpose, we restrict the admissible class to
supports containing a prescribed stabilizing core \(\omega\subset\Omega\).
We start with a result that relates the Riccati operators $\Pi_A,\Pi_B$ for $ A \subset B$.

\begin{lemma}
\label{lem:value_monotonicity_actuator}
Let \(A,B\subset\Omega\) be measurable actuator supports with \(A\subset B\)
up to a null set. Then, for every \(y_0\in H\), $V_B(y_0)\leq V_A(y_0)$. If the corresponding Riccati operators exist, then  $0\leq \Pi_B\leq \Pi_A$ in the sense of quadratic forms on $H$.
\end{lemma}

\begin{proof}
Let \(u\in L^2(0,\infty;L^2(\Omega))\) be any admissible control for the
support \(A\), and set \(v:=\chi_Au\). Since \(A\subset B\), we have
$\chi_Bv=\chi_B\chi_Au=\chi_Au$. Thus the support \(B\), used with the control \(v\), produces the same state trajectory as the support \(A\), used with the control \(u\). Moreover, \(\norm{v(t)}_{L^2(\Omega)}\leq\norm{u(t)}_{L^2(\Omega)}\) for a.e.\@ \(t\).
Hence, the cost for \(B\) is not larger than that for \(A\). Taking the
infimum over all \(u\) gives \(V_B(y_0)\leq V_A(y_0)\). The operator
inequality follows from the Riccati representation
\(V_A(y)=\ip{\Pi_Ay}{y}_H\) and \(V_B(y)=\ip{\Pi_By}{y}_H\).
\end{proof}

The next lemma records a quantitative version of the classical
Datko--Pazy criterion; see \cite[Thm.~V.1.8, p.~300]{EngelNagel2000}.
We include the proof to make the dependence on the two uniform constants
explicit.
\begin{lemma}
\label{lem:uniform_datko_estimate}
Let \(\{(T_\lambda(t))_{t\geq0}:\lambda\in\Lambda\}\) be a family of
strongly continuous semigroups on a Hilbert space \(H\). Suppose that
there are constants \(C\geq1\) and \(K\geq1\), independent of \(\lambda\),
such that
\[
    \int_0^\infty\|T_\lambda(t)z\|_H^2\dd t\leq C\|z\|_H^2
    \quad\lambda\in\Lambda,\ z\in H,
    \qquad
    \sup_{\lambda\in\Lambda}\sup_{0\leq t\leq1}
    \|T_\lambda(t)\|_{\Ls(H)}\leq K.
\]
Then, we have the following uniform exponential stability
\[
    \|T_\lambda(t)\|_{\Ls(H)}
    \leq K\sqrt C\,e^{1/(2C)}e^{-t/(2C)}
    \qquad \lambda\in\Lambda,\ t\geq0.
\]
\end{lemma}
\begin{proof}
For fixed \(\lambda\in\Lambda\) and \(z\in H\),  we define the remaining energy  $F(t):=\int_t^\infty\|T_\lambda(s)z\|_H^2\dd s$. Using the semigroup property, applied with initial state \(T_\lambda(t)z\), we can write
\begin{equation}
\label{eq:lem_expo}
        F(t)= \int_0^\infty\|T_\lambda(t+r)z\|_H^2\dd r = \int_0^\infty\|T_\lambda(r)T_\lambda(t)z\|_H^2\dd r
    \leq C\|T_\lambda(t)z\|_H^2=-CF'(t).
\end{equation}
Thus, we have  \(F'(t)\leq-F(t)/C\), and we can write
\[
    F(t)\leq F(0)e^{-t/C}\leq C e^{-t/C}\|z\|_H^2.
\]
Now, let \(t\geq1\)
and \(s\in[t-1,t]\). Since \(t-s\in[0,1]\), we obtain
\[
    \|T_\lambda(t)z\|_H^2
    =\|T_\lambda(t-s)T_\lambda(s)z\|_H^2
    \leq K^2\|T_\lambda(s)z\|_H^2.
\]
Integrating in \(s\) over $(t-1,t)$ and using \eqref{eq:lem_expo} yields the following
\[
    \|T_\lambda(t)z\|_H^2
    \leq K^2\int_{t-1}^t\|T_\lambda(s)z\|_H^2\dd s
    \leq K^2F(t-1)
    \leq K^2 C e^{-(t-1)/C}\|z\|_H^2.
\]
Taking square roots and the supremum over \(\|z\|_H=1\) proves the
estimate for \(t\geq1\). For \(0\leq t\leq1\), the bound
\(\|T_\lambda(t)\|\leq K\) gives the same estimate because
\(\sqrt C\,e^{(1-t)/(2C)}\geq1\). All constants are independent of
\(\lambda\).
\end{proof}

We are now in the position to prove that the assumptions of \cref{thm_sufficient_condition_for_under_uniform_remainder}
are satisfied for the admissible set $\Aadomega$.

\begin{proposition}
\label{cor:fixed_core_uniform_bounds}
Let \(\omega\subset\Omega\) be a fixed set such that $(\A,\chi_\omega)$ is exponentially stabilizable.
Assume \((Q,\A)\) is exponentially detectable.
Set
\[
     \Aadomega:=
    \{A \subset \Omega:\ \omega\subset A\text{ up to a null set}\},
\]
Then for every $A \in \Aadomega$, the pair $(\A,\chi_A)$ is exponentially stabilizable, i.e.,
$\Aadomega$ satisfies \cref{ass:lqr}. In addition, for all \(y_0\in H\), we have
\(V_A(y_0)\leq V_\omega(y_0)\leq \|\Pi_\omega\|_{\Ls(H)}\norm{y_0}_H^2\), and thus  \(0\leq\Pi_A\leq\Pi_\omega\) with \(\norm{\Pi_A}_{\Ls(H)}\leq\norm{\Pi_\omega}_{\Ls(H)}\).
%In addition, we have
%\[
%    V_A(y_0)\leq V_\omega(y_0)
%    \leq C_\omega\norm{y_0}_H^2,
%    \qquad C_\omega:=\|\Pi_\omega\|_{\Ls(H)},\quad y_0\in H,
%\]
%and
%\[
%    0\leq\Pi_A\leq\Pi_\omega,
%    \qquad
%    \norm{\Pi_A}_{\Ls(H)}\leq\norm{\Pi_\omega}_{\Ls(H)} .
%\]
Moreover, any exponentially stabilizing feedback supported
on \(\omega\) is an admissible stabilizing feedback for every
\(A\in\Aadomega\),
with the same closed-loop semigroup.
Thus, exponential stabilizability can be
chosen uniformly on \(\Aadomega\).
% Alternatively, if the uncontrolled semigroup satisfies
% \(\norm{S(t)}_{\Ls(H)}\leq Me^{-\gamma t}\), \(t\geq0\), with
% \(\gamma>0\), then the uniform value bound
% \[
%     V_A(y_0)
%     \leq
%     \frac{M^2}{2\gamma}\norm{y_0}_H^2
% \]
% holds for every actuator support \(A\), by choosing the zero control.
% In either case, the optimal closed-loop semigroups satisfy
In particular, there are \(M_{\cl},\gamma_{\cl}>0\) such that
\[
    \norm{S_A^{\cl}(t)}_{\Ls(H)}
    \leq M_{\cl}e^{-\gamma_{\cl}t},\qquad t\geq0, \ A \in \Aadomega.
\]
% with \(M_{\cl},\gamma_{\cl}>0\) independent of \(A\) in the respective class.
%
In addition, if $\Pi_\omega$ can be extended to an operator in $\Ls(L^1(\Omega),L^\infty(\Omega))$, then
$\Pi_A \in \Ls(L^1(\Omega),L^\infty(\Omega))$
with $\|\Pi_A\|_{\Ls(L^1(\Omega),L^\infty(\Omega))} \le\|\Pi_\omega\|_{\Ls(L^1(\Omega),L^\infty(\Omega))} $
for all $A \in \Aadomega$.
\end{proposition}

\begin{proof}
Let $K_\omega:H\to L^2(\Omega)$ be the optimal feedback associated to $\omega$, which is exponentially stabilizing.
For $A\in\Aadomega$ we have
\(\A+\chi_AK_\omega=\A+\chi_\omega K_\omega\), and the same exponentially stable
closed-loop semigroup stabilizes all larger actuator supports. Thus, the
standard LQR theorem (\cref{thm:riccati}) then gives well-posedness for every
\(A\in\Aadomega\).
Then \cref{lem:value_monotonicity_actuator} gives \(V_A(y_0)\leq V_\omega(y_0)\).

The Riccati bounds follow again from
the quadratic value representation in \cref{lem:value_monotonicity_actuator} and positivity of the Riccati operators.
% In
% the exponentially stable uncontrolled case, the admissible choice \(u=0\)
% gives
% \[
%     V_A(y_0)
%     \leq
%     \int_0^\infty \norm{S(t)y_0}_H^2\dd t
%     \leq
%     \frac{M^2}{2\gamma}\norm{y_0}_H^2 .
% \]
To obtain uniform stability of the optimal feedbacks, we verify the two
hypotheses of \cref{lem:uniform_datko_estimate}.

By exponential detectability, we can use the same observer $L$ for every support \(A\).
The operator $\A+LQ$ generates an exponentially stable semigroup. For $z\in H$, let $y(t):=S_A^{\cl}(t)z$ with corresponding optimal feedback control $u_A$.
We write $\dot y = (\A+LQ)y -LQ y + \chi_Au_A$. Arguing as in the proof of \cite[Thm. 6.2.7]{CurtainZwart1995}, we bound $\norm{y}_{L^2(0,\infty;H)}$ in terms of $\norm{z}_H$,
$\norm{Qy}_{L^2(0,\infty;H_Q)}$, and $\norm{\chi_Au_A}_{L^2(0,\infty;H)}$. The last two terms are controlled by the optimal value, yielding
% Duhamel's formula and
% Young's inequality give
\[
    \int_0^\infty\|S_A^{\cl}(t)z\|_H^2\,\dd t
    \leq  C\bigl(\|z\|_H^2+V_A(z)\bigr)
    \leq C_\omega\|z\|_H^2,
\]
where \(C\) is independent of $A$, % only on the fixed problem data
and the last inequality
follows from \(V_A(z)\leq V_\omega(z)\), with \(C_\omega\geq1\). This
verifies the first hypothesis of \cref{lem:uniform_datko_estimate}
uniformly in \(A\).

Furthermore, the closed-loop generator \(\A_A^{\cl}=\A+\KA\), where
\(\KA=-\alpha^{-1}\chi_A\Pi_A\), is a bounded perturbation satisfying
\[
    \|\KA\|_{\Ls(H)}
    \leq\alpha^{-1}\|\chi_A\Pi_A\|_{\Ls(H)}
    \leq  \alpha^{-1} \|\Pi_\omega\|_{\Ls(H)}.
\]
% \end{samepage}
Choose \(M_0\geq1\) and \(\eta\in\R\) such that the uncontrolled
semigroup, i.e., the semigroup generated by $\A$, satisfies \(\|S(t)\|\leq M_0e^{\eta t}\) for \(t\geq0\).
The bounded perturbation theorem
\cite[Ch.~3, Thm.~1.1, p.~76]{Pazy1983} then gives
\[
    \|S_A^{\cl}(t)\|_{\Ls(H)}
    \leq M_0 e^{(\eta+M_0\|\KA\|_{\Ls(H)})t}
    \leq M_0 e^{(\eta+M_0\alpha^{-1}\|\Pi_\omega\|_{\Ls(H)})t}.
\]
Consequently, the common constant $K:=M_0\exp\left(\max\{0,\eta+M_0\alpha^{-1}\|\Pi_\omega\|_{\Ls(H)}\}\right)\geq 1$ satisfies
\(\sup_A\sup_{0\leq t\leq1}\|S_A^{\cl}(t)\|_{\Ls(H)}\leq K\). Applying  \cref{lem:uniform_datko_estimate} to the family
\( (S_A^{\cl})_{A\in \Aadomega}\) now proves the asserted exponential decay with $\gamma_{\cl}=\frac1{2C_\omega}$ and $M_{\cl}=K\sqrt{C_\omega}\,e^{1/(2C_\omega)}$
independent of $A \in \Aadomega$.

Since $\Pi_A$ is self-adjoint and nonnegative for $A \in \Aadomega$, we have for $v,w\in H$,
\begin{multline*}
 ( \Pi_Av,w)^2 \le (\Pi_Av,v)(\Pi_Aw,w) \le (\Pi_\omega v,v)(\Pi_\omega w,w) \\
 \le \|\Pi_\omega\|_{\Ls(L^1(\Omega),L^\infty(\Omega))}^2 \|v\|_{L^1(\Omega)}^2 \|w\|_{L^1(\Omega)}^2.
\end{multline*}
As $L^2(\Omega)$ is dense in $L^1(\Omega)$ it follows that $\|\Pi_A\|_{\Ls(L^1(\Omega),L^\infty(\Omega))} \le\|\Pi_\omega\|_{\Ls(L^1(\Omega),L^\infty(\Omega))} $.
\end{proof}

For the control of the parabolic partial differential equation (cf., \cref{prop:parabolic_is_nice}),
we have the following result on the regularity of the Riccati operators.
It implies that the operators $\Pi_A$ are uniformly bounded in $\Ls(H,L^\infty(\Omega))$ for $A\in \Aadomega$ with $\Aadomega$ as in \cref{cor:fixed_core_uniform_bounds}.

\begin{lemma}
\label{lem:Riccati_linfty_estimate}
 Let the control problem be as in \cref{prop:parabolic_is_nice}.
 Then there is $C>0$ such that for every measurable support
 \(A\subset\Omega\) for which \((\A,\chi_A)\) is exponentially stabilizable,
 we have $\Pi_A \in \Ls(H,L^\infty(\Omega))$ and
 \[
  \norm{\Pi_A}_{\Ls(H,L^\infty(\Omega))} \le  C(1+ \norm{\Pi_A}_{\Ls(H)}^2).
 \]
\end{lemma}
\begin{proof}
Let \(S(t)=e^{t\A}\) denote the semigroup on \(H=L^2(\Omega)\)
generated by the uncontrolled realization of
\(\A=\kappa\Delta-b\cdot\nabla-c\) with the fixed boundary condition
\(\mathcal G y=0\) from \eqref{eq:motivating_pde}. Its adjoint \(S(t)^*\) is generated by \(\A^*\), understood through the adjoint form with its associated boundary condition. The elliptic form and its adjoint have the same weak coercivity bound.
Hence, by \cite[Thm.~6.8]{Ouhabaz2005} or
\cite[Thm.~5.2 and Rem.~5.4, pp.~30--31]{Daners2000}, there exists
\(C_{\mathrm{sm}}>0\) such that
\[
    \norm{S(t)^*}_{\Ls(H,L^\infty(\Omega))}
    \leq C_{\mathrm{sm}}t^{-d/4},\qquad 0<t\leq1,
\]
for both boundary conditions, where \(C_{\mathrm{sm}}\) depends only on the fixed PDE data.

The operator \(Q_A:= Q^*Q -\alpha^{-1}\Pi_A\chi_A\Pi_A\) satisfies \(\norm{Q_A}_{\Ls(H)}\leq \norm{Q}_{\Ls(H,H_Q)}^2+\alpha^{-1}\norm{\Pi_A}_{\Ls(H)}^2\).
Integrating the weak Riccati equation along \(S(t)\) over \([0,1]\), we obtain the identity on \(\Ls(H)\)
\begin{equation} \label{eq_integrated_riccati_Linfty}
    \Pi_A=S(1)^*\Pi_A S(1)
    +\int_0^1 S(t)^*Q_A S(t)\dd t.
\end{equation}
Define \(M_1:=\sup_{0\leq t\leq1}\norm{S(t)}_{\Ls(H)}\). Then we can estimate
\[
\begin{aligned}
    \norm{S(1)^*\Pi_A S(1)}_{\Ls(H,L^\infty(\Omega))}
    &\leq C_{\mathrm{sm}}\norm{\Pi_A S(1)}_{\Ls(H)}
    \leq C_{\mathrm{sm}}\norm{\Pi_A}_{\Ls(H)}M_1 ,\\
    \norm{S(t)^*Q_A S(t)}_{\Ls(H,L^\infty(\Omega))}
    &\leq C_{\mathrm{sm}}t^{-d/4}\norm{Q_A S(t)}_{\Ls(H)} \\
    &\leq C_{\mathrm{sm}}M_1
    ( \norm{Q}_{\Ls(H,H_Q)}^2+\alpha^{-1}\norm{\Pi_A}_{\Ls(H)}^2)t^{-d/4}.
\end{aligned}
\]
Since \(d\leq3\), the last bound is integrable in time and shows that the \(\Ls(H)\)-valued integral in \eqref{eq_integrated_riccati_Linfty}
also belongs to \(\Ls(H,L^\infty(\Omega))\). The identity \eqref{eq_integrated_riccati_Linfty} therefore gives directly
\[
\begin{aligned}
    \norm{\Pi_A }_{\Ls(H,L^\infty(\Omega))}
    &\leq\norm{S(1)^*\Pi_A S(1)}_{\Ls(H,L^\infty(\Omega))}
    +\int_0^1\norm{S(t)^*Q_A S(t)}_{\Ls(H,L^\infty(\Omega))}\dd t\\
    &\leq C_{\mathrm{sm}}M_1
    \left(\norm{\Pi_A}_{\Ls(H)}+
    \frac{\norm{Q_A}_{\Ls(H)}}{1-d/4}\right)
    \le C(1+ \norm{\Pi_A}_{\Ls(H)}^2),
\end{aligned}
\]
where $C$ is independent of $A$.
\end{proof}

We now consider \cref{alg:lyapunov_support_descent} with
\(\Aad=\Aadomega\) as in \cref{cor:fixed_core_uniform_bounds}. Every
admissible support contains the prescribed stabilizing core \(\omega\) up
to a null set. Accordingly, \cref{thm:topological_derivative,thm:pmp} apply
on \(\Omega\setminus\omega\), with \(A\setminus(B_\varepsilon(x)\setminus\omega)\)
as the removal trial.

The switching function \(G_A\) is unchanged. In \eqref{eq:stationarity_defect},
we use \(q_A=\chi_{\Omega\setminus\omega}\min\{\Theta_A,0\}\) and keep the
definitions of \(\eta_A\) and \(R_A\). Since
\(D_{k,s}\subset R_k\subset\Omega\setminus\omega\), every trial preserves
the core and is admissible. The other algorithmic steps are unchanged.
In \cref{lem:computable_error_estimate_V,thm_convergence_under_uniform_remainder,thm_sufficient_condition_for_under_uniform_remainder},
all comparisons, remainder bounds, and infima are taken over \(\Aadomega\).
Without a prescribed core, the original algorithm is recovered.

\begin{theorem}
Let the control problem be as in \cref{prop:parabolic_is_nice}.
Suppose $\Aad =  \Aadomega$ with $ \Aadomega$ as in \cref{cor:fixed_core_uniform_bounds}.
Let $(A_k)$ be an infinite sequence of iterates of \cref{alg:lyapunov_support_descent}
with the core-preserving updates described above.
Then the assumptions of \cref{thm_convergence_under_uniform_remainder,thm_sufficient_condition_for_under_uniform_remainder}
for this restricted class are satisfied. In particular, $(A_k)$ is a minimizing sequence of $J$ defined in \eqref{eq:motivating_support_cost}
over \(\Aadomega\).
\end{theorem}
\begin{proof}
\Cref{cor:fixed_core_uniform_bounds} shows that the first condition of \cref{thm_sufficient_condition_for_under_uniform_remainder} is satisfied.
The uniform bound on the $\Ls(H)$-norm of the Riccati operators from \cref{cor:fixed_core_uniform_bounds} together with the $\Ls(H,L^\infty(\Omega))$ estimate of \cref{lem:Riccati_linfty_estimate}
shows that the second condition of  \cref{thm_sufficient_condition_for_under_uniform_remainder} is satisfied.
Thus the uniform remainder estimate assumption of \cref{thm_convergence_under_uniform_remainder} is verified, and the convergence statement follows.
\end{proof}

\section{Numerical experiments}
 \label{sec:numerics}
We apply  \cref{alg:lyapunov_support_descent}
to optimize actuator location and shape for an unstable parabolic system.
Starting from a common support $A_0$, the algorithm reduces the combined
LQR and support cost, with lower support weights assigned to a preferred
cross, annulus, disk, or with a uniform support weight throughout the domain. We examine the resulting shapes, the objective
and stationarity histories, and the closed-loop response.

\begin{figure}[!htb]
 \centering
 \begin{minipage}[c]{0.4\linewidth}
  \includegraphics[width=1.\linewidth,height=3.5cm,
  keepaspectratio]{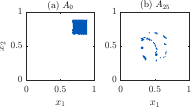}
 \end{minipage}\hspace{0.03\linewidth}%
 \begin{minipage}[c]{0.39\linewidth}
  \caption{Actuator supports: (a) common initial support;
  (b) final support for \(\beta\equiv2\).}
  \label{fig:numerical_initial_support}
  \label{fig:numerical_uniform_support}
 \end{minipage}
\end{figure}

\paragraph{Model and discretization}
We consider \eqref{eq:motivating_pde} on \(\Omega=(0,1)^2\) with
homogeneous Dirichlet conditions, \(\kappa=0.1\), \(\alpha=10^{-3}\), $c(x):=-2-\sin(\pi x_1)\sin(\pi x_2)$, and
\[
\begin{aligned}
% c(x)&:=-2-\sin(\pi x_1)\sin(\pi x_2),\\
 b(x)&:=\bigl(0.30+0.15\sin(2\pi x_2),\
             0.10+0.15\cos(2\pi x_1)\bigr),\\
 y_0(x)&:=\sin(\pi x_1)\sin(\pi x_2)
          +0.25\sin(2\pi x_1)\sin(\pi x_2).
\end{aligned}
\]
The uncontrolled system is exponentially unstable, whereas
the LQR feedback for the common initial support
\(A_0=[0.68,0.90]^2\) is exponentially stabilizing
(\cref{fig:numerical_trajectories}). This support is sampled at triangle barycentres and shown in \cref{fig:numerical_initial_support}(a).

For the spatial discretization, we use continuous piecewise-linear finite elements for the state and piecewise-constant functions on the triangles for the control. The triangulation \(\mathcal T_h\) is generated with
the target maximum edge length \(h_{\max}=0.02\) and contains \(2774\)
free state nodes and \(5746\) triangles.
Supports are unions of triangles; the \(L^2\)-costs use a consistent state
mass matrix and an \(L^2\)-orthonormal control basis.
All computations are implemented in MATLAB. The routines \texttt{icare}
and \texttt{lyap} solve the generalized algebraic Riccati and Lyapunov
equations used to evaluate the discrete infinite-horizon LQR value
\(V_h(A)\) and switching density. The routine \texttt{icare} also returns
the discrete closed-loop eigenvalues \(\lambda_j\). Their largest real
part, \(s_h^{\rm cl}:=\max_j\operatorname{Re}\lambda_j\), is the spectral
abscissa; \(s_h^{\rm cl}<0\) characterizes exponential stability of the
discrete closed-loop system. We denote the combined objective by
\(J_h(A)\) and display trajectories for \(t\in[0,8]\).

For each triangle \(K\in\mathcal T_h\), let \(\Theta_{k,K}\) denote the
discrete flip indicator at iteration \(k\), and set
\(q_{k,K}:=\min\{\Theta_{k,K},0\}\).
We rank triangles with \(q_{k,K}<0\) from the most negative value upwards
and backtrack over prefixes of this list, starting with all violating
triangles. The area reduction factor is \(\rho = 1/2\), the Armijo parameter is
\( \sigma = 10^{-4} \), and at most \(40\) backtracking reductions are allowed per
search. With \(|K|\) denoting triangle area, we record
\[
 \eta_{k,h}=\sum_{K\in\mathcal T_h}|K|\,|q_{k,K}|
          =-\sum_{K\in\mathcal T_h}|K|q_{k,K},
 \qquad
 \eta_{\infty,k,h}=\max_{K\in\mathcal T_h}|q_{k,K}|.
\]
At each iteration, we check the stationarity test
\(\eta_{\infty,k,h}\leq10^{-4} \eqqcolon \varepsilon_{\rm stop}  \), instead of the integral test in
\cref{alg:lyapunov_support_descent}. A run terminates if this test
is satisfied, if no tested ranked-prefix update is accepted, or after
\(25\) updates.

\paragraph{Actuator designs}

\begin{figure}[!htb]
 \centering
 \includegraphics[width=0.66\linewidth,height=5.5cm,
  keepaspectratio]{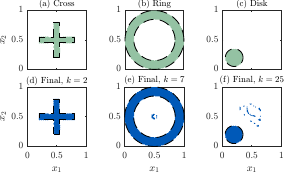}
 \caption{Preferred regions (green, top) and final actuator supports
 (blue, bottom).}
 \label{fig:numerical_supports}
\end{figure}

For each preferred region \(D\), we set \(\beta=10^{-7}\) on \(D\) and
\(\beta=2\) outside, sampled at triangle barycentres. This weighting
encourages placement within \(D\), while allowing portions outside \(D\)
when beneficial for the combined objective.
We take \(\omega=\varnothing\) and start all runs from \(A_0\). \Cref{fig:numerical_supports} shows the preferred regions in (a)--(c) and the final supports in (d)--(f). The black dashed curves mark the exact boundaries of the preferred regions;
the green regions are their triangular-mesh
approximations. Some final supports retain small isolated components,
which are included in both the plots and the reported costs.
For comparison, we set \(\beta\equiv2\) on \(\Omega\), which corresponds to
\(D=\varnothing\).
The final support consists of small disconnected components in the
interior of the domain (\cref{fig:numerical_uniform_support}(b)).
Its total area is smaller than in the three spatially weighted cases
(see \cref{tab:numerical_costs}), reflecting the penalty on the entire
actuator support.

\subsection{Objective and stationarity}
\begin{figure}[!htbp]
 \centering
 \includegraphics[width=0.92\linewidth,height=4.3cm,
  keepaspectratio]{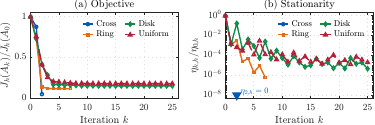}
 \caption{Normalized objectives and relative stationarity defects.}
 \label{fig:numerical_histories}
\end{figure}

Let \(k_{\rm f}\) denote the final iteration index and \(A_{k_{\rm f}}\)
the resulting support. \Cref{tab:numerical_costs} reports
the costs and support areas.
\begin{table}[htbp]
 \centering\small
 \setlength{\tabcolsep}{4pt}
 \begin{tabular}{lrrrrrr}\toprule
 Region $D$ & $k_{\rm f}$ & $J_h(A_0)$ & $J_h(A_{k_{\rm f}})$ & $V_h(A_{k_{\rm f}})$ & $\int_{A_{k_{\rm f}}}\beta$ & $|A_{k_{\rm f}}|$ \\
 \midrule
 Cross & 2 & 0.42861 & 0.02129 & 0.02129 & \num{1.115e-08} & 0.11151 \\
 Ring & 7 & 0.35823 & 0.04175 & 0.03459 & \num{7.152e-03} & 0.40640 \\
 Off-centre disk & 25 & 0.42861 & 0.06527 & 0.04624 & \num{1.903e-02} & 0.08004 \\
 \(\varnothing\) (uniform) & 25 & 0.42861 & 0.07679 & 0.05091 & \num{2.588e-02} & 0.01294 \\
 \bottomrule\end{tabular}
 \caption{Numerical results for the four actuator designs.}
 \label{tab:numerical_costs}
\end{table}
\Needspace{5\baselineskip}
For the disk, iteration two already yields a support entirely within
\(D\). Further optimization reduces the objective, although the final
support retains additional components outside \(D\).
\Cref{fig:numerical_histories} shows the monotone decrease in the objectives.

The cross satisfies the discrete stationarity test after two updates:
\(\Theta_{2,K}>0\) on every triangle, hence \(\eta_{2,h}=0\),
indicated by the downward marker in \cref{fig:numerical_histories}.
For the ring, disk, and uniform penalty, \(\eta_{k,h}\) decreases
substantially, but the stationarity test is not satisfied. The ring stops
after seven updates when the line search accepts no tested update,
whereas the disk and uniform-penalty runs reach the limit of \(25\) updates. \Cref{tab:numerical_stationarity} reports the residuals and termination criteria.
The elapsed optimization times were approximately
\(1.77\), \(13.97\), \(62.07\), and \(60.02\) hours for the cross,
ring, disk, and uniform-penalty examples, respectively.

\begin{table}[htbp]
 \centering\small
 \setlength{\tabcolsep}{3.25pt}
 \begin{tabular}{l *{3}{S[table-format=1.3e-1]} S[table-format=-1.3] l}\toprule
 Region $D$ & {$\eta_{k_{\rm f},h}$} & {$\eta_{k_{\rm f},h}/\eta_{0,h}$} & {$\eta_{\infty,k_{\rm f},h}$} & {$s_h^{\rm cl}$} & Termination \\
 \midrule
 Cross & 0.000 & 0.000 & 0.000 & -6.481 & Stationarity \\
 Ring & 1.320e-04 & 6.201e-07 & 2.878e-02 & -7.234 & No accepted step \\
 Off-centre disk & 9.942e-04 & 4.680e-06 & 1.029e-01 & -5.095 & Iteration limit \\
 \(\varnothing\) (uniform) & 4.242e-03 & 1.998e-05 & 1.600e-01 & -5.068 & Iteration limit \\
 \bottomrule\end{tabular}
 \caption{Final stationarity diagnostics and closed-loop spectral
 abscissae \(s_h^{\rm cl}\).}
 \label{tab:numerical_stationarity}
\end{table}

\begin{figure}[!htbp]
 \centering
 \includegraphics[width=0.92\linewidth,height=4.5cm,
  keepaspectratio]{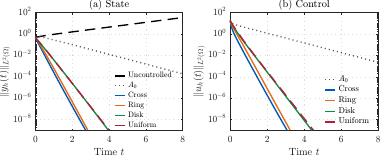}
 \caption{Evolution of the \(L^2\)-norms of the state and control.}
 \label{fig:numerical_trajectories}
\end{figure}

\paragraph{Closed-loop behavior}
\Cref{fig:numerical_trajectories} shows the evolution of the
\(L^2\)-norms of the state and control for the common initial support
\(A_0\) and the final supports in all four examples.
For the chosen initial state, the final-support feedback controls yield faster
decay than the feedback for \(A_0\), whereas the uncontrolled state grows.
The spectral abscissae are negative for \(A_0\) and all final supports,
confirming exponential stability of the discrete closed-loop systems;
\Cref{tab:numerical_stationarity} presents the final values.
The disk and uniform-penalty cases exhibit similar decay of
the state and control norms, consistent with their close closed-loop
spectral abscissae. The relative Riccati and Lyapunov residuals at the final iterations are below  \(5\times10^{-12}\) and \(8\times10^{-13}\), respectively, supporting the accuracy of the algebraic solves.

\begin{figure}[!htbp]
 \centering
 \setlength{\abovecaptionskip}{6pt}
 \includegraphics[width=\linewidth,height=4cm,
  keepaspectratio]{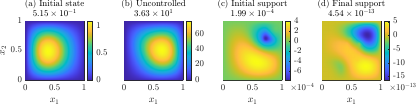}\par
 \nointerlineskip
 \caption{State snapshots for the disk example: (a) \(y_0\) at \(t=0\);
 (b)--(d) states at \(t=8\) without control, with the initial-support feedback, and with the
 final-support feedback. Separate color bars show each panel's scale;
 the second title lines give the \(L^2\)-norms.}
 \label{fig:disk_state_snapshots}
\end{figure}

\Cref{fig:disk_state_snapshots} compares the states for the disk example at \(t=0\)
and \(t=8\), illustrating the growth of the uncontrolled state and the stronger stabilization achieved with the final support.

%  \section*{Conclusion} We developed a topological method for actuator-support optimization in infinite-horizon LQ parabolic control. The Riccati-based sensitivity formula and its Lyapunov representation lead to a measure-controlled Armijo algorithm. Under the stated uniform remainder and admissibility assumptions, the admissible first-order stationarity defect tends to zero. Numerical experiments illustrate actuator relocation, changes in topology, and objective reduction while preserving closed-loop stability.

\section*{Acknowledgments}
After completing the initial draft, the authors used OpenAI Codex to proofread the manuscript, check the citations, and assist in reviewing the mathematical arguments. The authors assume responsibility for all content, including ideas, mathematical results and proofs, algorithms, and numerical experiments.

\AddToHookNext{env/thebibliography/begin}{\setlength{\itemsep}{0pt}}
\bibliographystyle{siamplain}
\Needspace{6\baselineskip}
\bibliography{Optimal_placment}

@book{LasieckaTriggiani2000,
  author    = {Lasiecka, Irena and Triggiani, Roberto},
  title     = {Control Theory for Partial Differential Equations: Continuous and Approximation Theories. I. Abstract Parabolic Systems},
  series    = {Encyclopedia of Mathematics and its Applications},
  volume    = {74},
  publisher = {Cambridge University Press},
  year      = {2000}
}

@book{Pazy1983,
  author    = {Pazy, Amnon},
  title     = {Semigroups of Linear Operators and Applications to Partial Differential Equations},
  series    = {Applied Mathematical Sciences},
  volume    = {44},
  publisher = {Springer},
  year      = {1983},
  doi       = {10.1007/978-1-4612-5561-1}
}

@book{Lunardi1995,
  author    = {Lunardi, Alessandra},
  title     = {Analytic Semigroups and Optimal Regularity in Parabolic Problems},
  series    = {Progress in Nonlinear Differential Equations and Their Applications},
  volume    = {16},
  publisher = {Birkh{\"a}user},
  year      = {1995}
}

@book{EngelNagel2000,
  author    = {Engel, Klaus-Jochen and Nagel, Rainer},
  title     = {One-Parameter Semigroups for Linear Evolution Equations},
  series    = {Graduate Texts in Mathematics},
  volume    = {194},
  publisher = {Springer},
  year      = {2000}
}

@book{LionsMagenes1972,
  author    = {Lions, Jacques-Louis and Magenes, Enrico},
  title     = {Non-Homogeneous Boundary Value Problems and Applications},
  publisher = {Springer},
  year      = {1972}
}

@book{BensoussanDaPratoDelfourMitter2007,
  author    = {Bensoussan, Alain and Da Prato, Giuseppe and Delfour, Michel C. and Mitter, Sanjoy K.},
  title     = {Representation and Control of Infinite Dimensional Systems},
  edition   = {2},
  series    = {Systems \& Control: Foundations \& Applications},
  publisher = {Birkh{\"a}user},
  address   = {Boston},
  year      = {2007},
  doi       = {10.1007/978-0-8176-4581-6}
}

@book{TucsnakWeiss2009,
  author    = {Tucsnak, Marius and Weiss, George},
  title     = {Observation and Control for Operator Semigroups},
  series    = {Birkh{\"a}user Advanced Texts},
  publisher = {Birkh{\"a}user},
  address   = {Basel},
  year      = {2009},
  doi       = {10.1007/978-3-7643-8994-9}
}

@article{Wachsmuth2024,
  author  = {Wachsmuth, Daniel},
  title   = {A topological derivative-based algorithm to solve optimal control problems with {$L^0(\Omega)$} control cost},
  journal = {Journal of Nonsmooth Analysis and Optimization},
  volume  = {5},
  year    = {2024},
  pages   = {12366},
  doi     = {10.46298/jnsao-2024-12366}
}

@book{CurtainZwart1995,
  author    = {Curtain, Ruth F. and Zwart, Hans},
  title     = {An Introduction to Infinite-Dimensional Linear Systems Theory},
  series    = {Texts in Applied Mathematics},
  volume    = {21},
  publisher = {Springer},
  year      = {1995},
  doi       = {10.1007/978-1-4612-4224-6}
}

@book{Zabczyk1992,
  author    = {Zabczyk, Jerzy},
  title     = {Mathematical Control Theory: An Introduction},
  publisher = {Birkh{\"a}user},
  address   = {Boston},
  year      = {1992}
}

@book{Ouhabaz2005,
  author    = {Ouhabaz, El Maati},
  title     = {Analysis of Heat Equations on Domains},
  series    = {London Mathematical Society Monographs Series},
  volume    = {31},
  publisher = {Princeton University Press},
  address   = {Princeton, NJ},
  year      = {2005}
}

@article{FernandezCaraGuerrero2006,
  author  = {Fernandez-Cara, Enrique and Guerrero, Sergio},
  title   = {Global {Carleman} inequalities for parabolic systems and applications to controllability},
  journal = {SIAM Journal on Control and Optimization},
  volume  = {45},
  number  = {4},
  pages   = {1395--1446},
  year    = {2006},
  doi     = {10.1137/S0363012904439696}
}

@article{FernandezCaraGonzalezBurgosGuerreroPuel2006,
  author  = {Fernandez-Cara, Enrique and Gonzalez-Burgos, Manuel and Guerrero, Sergio and Puel, Jean-Pierre},
  title   = {Null controllability of the heat equation with boundary {Fourier} conditions: the linear case},
  journal = {ESAIM: Control, Optimisation and Calculus of Variations},
  volume  = {12},
  number  = {3},
  pages   = {442--465},
  year    = {2006},
  doi     = {10.1051/cocv:2006010}
}

@book{FursikovImanuvilov1996,
  author    = {Fursikov, A. V. and Imanuvilov, O. Yu.},
  title     = {Controllability of Evolution Equations},
  series    = {Lecture Notes Series},
  volume    = {34},
  publisher = {Seoul National University, Research Institute of Mathematics, Global Analysis Research Center},
  address   = {Seoul},
  year      = {1996}
}

@article{Morris2011,
  author  = {Morris, Kirsten A.},
  title   = {Linear-Quadratic Optimal Actuator Location},
  journal = {IEEE Transactions on Automatic Control},
  volume  = {56},
  number  = {1},
  pages   = {113--124},
  year    = {2011},
  doi     = {10.1109/TAC.2010.2052151}
}

@article{DarivandiMorrisKhajepour2013,
  author  = {Darivandi, Neda and Morris, Kirsten A. and Khajepour, Amir},
  title   = {An Algorithm for {LQ} Optimal Actuator Location},
  journal = {Smart Materials and Structures},
  volume  = {22},
  number  = {3},
  pages   = {035001},
  year    = {2013},
  doi     = {10.1088/0964-1726/22/3/035001}
}

@article{EdalatzadehMorris2019,
  author  = {Edalatzadeh, M. Sajjad and Morris, Kirsten A.},
  title   = {Optimal Actuator Design for Semilinear Systems},
  journal = {SIAM Journal on Control and Optimization},
  volume  = {57},
  number  = {4},
  pages   = {2992--3020},
  year    = {2019},
  doi     = {10.1137/18M1171229}
}

@article{KaliseKunischSturm2018,
  author  = {Kalise, Dante and Kunisch, Karl and Sturm, Kevin},
  title   = {Optimal Actuator Design Based on Shape Calculus},
  journal = {Mathematical Models and Methods in Applied Sciences},
  volume  = {28},
  number  = {13},
  pages   = {2667--2717},
  year    = {2018},
  doi     = {10.1142/S0218202518500586}
}

@article{EdalatzadehKaliseMorrisSturm2022,
  author  = {Edalatzadeh, M. Sajjad and Kalise, Dante and Morris, Kirsten A. and Sturm, Kevin},
  title   = {Optimal Actuator Design for the {Euler--Bernoulli} Vibration Model Based on {LQR} Performance and Shape Calculus},
  journal = {IEEE Control Systems Letters},
  volume  = {6},
  pages   = {1334--1339},
  year    = {2022},
  doi     = {10.1109/LCSYS.2021.3093215}
}

@article{Daners2000,
  author  = {Daners, Daniel},
  title   = {Heat Kernel Estimates for Operators with Boundary Conditions},
  journal = {Mathematische Nachrichten},
  volume  = {217},
  pages   = {13--41},
  year    = {2000},
  doi     = {10.1002/1522-2616(200009)217:1<13::AID-MANA13>3.0.CO;2-6}
}
% \endgroup

\end{document}